\documentclass[12pt,reqno]{amsart}
\usepackage{amsmath}
\usepackage{graphicx}

\def\squarebox#1{\hbox to #1{\hfill\vbox to #1{\vfill}}}

\newcommand{\Z}{{\mathbb Z}}

\newcommand{\R}{{\mathbb R}}
\newcommand{\C}{{\mathbb C}}
\newcommand{\N}{{\mathbb N}}

\newcommand{\ttt}{|\hspace{-0.25mm}|\hspace{-0.25mm}|}
\renewcommand{\Re}{\mathop{\rm Re}\nolimits}

\theoremstyle{plain}

\newtheorem{thm}{Theorem}

\newtheorem{cor}{Corollary}

\newtheorem{prop}{Proposition}

\newtheorem{lem}{Lemma}

\newtheorem{rem}{Remark}

\newtheorem{defn}{Definition}

\def\clip{C^{\mbox{\footnotesize \rm Lip}}}
\def\Lip{\mbox{\rm Lip}}
\def\lip{\mbox{\footnotesize\rm Lip}}
\def\Int{\mbox{\rm Int}}
\def\i{{ i}}

\def\W{{\mathcal W}}
\def\endofproof{{\rule{6pt}{6pt}}}
\def\di{\displaystyle}

\begin{document}

\title[Distribution of periods of closed trajectories]
{Distribution of periods of closed trajectories\\ in exponentially shrinking intervals}

\author[V. Petkov]{Vesselin Petkov}
\address{Universit\'e Bordeaux I, Institut de Math\'ematiques de Bordeaux, 351, 
Cours de la Lib\'eration, 33405  Talence, France}
\email{petkov@math.u-bordeaux1.fr}
\author[L. Stoyanov]{Luchezar Stoyanov}
\address{University of Western Australia, School of Mathematics and
Statistics,  Perth, WA 6009,  Australia}
\email{stoyanov@maths.uwa.edu.au}
\thanks{The first author was partially supported by the ANR project NONAA}

\def\ts{\tilde{\sigma}}
\def\ss{{\mathcal S}}
\def\aa{{\mathcal A}}
\def\R{{\mathbb R}}
\def\iu{\underline{i}}
\vspace*{0,8cm}
\def\saa{\Sigma_A^+}
\def\sAA{\Sigma_{\aa}^+}
\def\sA{\Sigma_{\aa}}
\def\saa{\Sigma_A^+}
\def\sa{\Sigma_A}
\def\san{\Sigma^-_A}

\def\Prf{\mbox{\footnotesize\rm Pr}}
\def\Pr{\mbox{\rm Pr}}
\def\lc{{\mathcal L}}
\def\pxi{\phi + (\xi + i u)\psi_p}
\def\sa{\Sigma_A}
\def\ssa{\Sigma_A^+}
\def\ssn{\sum_{\sigma^n x = x}}
\def\lc{{\mathcal L}}
\def\lo{{\mathcal O}}
\def\oo{{\mathcal O}}
\def\ff{{\mathcal F}}
\def\mt{\Lambda}
\def\ep{\epsilon}
\def\ms{\medskip}
\def\bs{\bigskip}
\def\diam{\mbox{\rm diam}}
\def\rr{\mathcal R}
\def\pp{\mathcal P}
\def\w{{\sf w}}
\def\hU{\widehat{U}}
\def\hz{\hat{z}}
\def\tz{\tilde{z}}
\def\ty{\tilde{y}}
\def\be{\begin{equation}}
\def\ee{\end{equation}}
\def\beqn{\begin{eqnarray}}
\def\eeqn{\end{eqnarray}}
\def\beqn*{\begin{eqnarray*}}
\def\eeqn*{\end{eqnarray*}}
\def\e{\emptyset}
\def\hR{\widehat{R}}

\def\dk{\partial K}
\def\mtb{\mt_{\dk}}
\def\tf{\tilde{f}}
\def\tg{\tilde{g}}
\def\tq{\tilde{q}}
\def\tde{\tilde{\delta}}
\def\tp{\tilde{p}}
\def\td{\tilde{d}}
\def\tP{\widetilde{P}}
\def\tI{\widetilde{I}}
\def\en{\epsilon_n}
\def\dist{\mbox{\rm dist}}

\def\H{{\mathbb H}}

\maketitle
\begin{abstract} For hyperbolic flows over basic sets  we study the asymptotic of the number of closed trajectories $\gamma$ with periods $T_{\gamma}$ lying in exponentially shrinking intervals 
$(x - e^{-\delta x}, x + e^{-\delta x}),\:\delta > 0,\: x \to + \infty.$ 
A general result is established which concerns hyperbolic flows admitting 
symbolic models whose corresponding Ruelle transfer operators 
satisfy some spectral estimates. This result applies to
a variety of hyperbolic flows on basic sets, in particular to geodesic flows on manifolds of 
constant negative curvature and to open billiard flows.

\end{abstract}
\section{Introduction}
\renewcommand{\theequation}{\arabic{section}.\arabic{equation}}
\setcounter{equation}{0}

The purpose of this paper is to study the asymptotic behavior of the number of closed trajectories for hyperbolic flows $\varphi_t$ in compact invariant sets.  It is known that if $\pi(x)$ is the number of closed orbits with primitive
period (length) not greater than $x$,  we have the asymptotic
$$\lim_{x \to +\infty}  \frac{1}{x} \log \pi(x) = h_T \, ,$$
where $h_T > 0$ is the topological entropy of the flow $\varphi_t.$ To get more precise results one has to impose some conditions 
on the flow. Thus, if $\varphi_t$ is a weak-mixing Axiom A flow restricted to a non-trivial basic set, Parry and Pollicott proved  
\cite{kn:PP} that
\begin{equation} \label{eq:1.1}
\pi(x) \sim \frac{e^{h_T x}}{h_T x} \quad , \quad x \to +\infty.
\end{equation}
This asymptotic generalizes the classical result of Margulis \cite{kn:M} for geodesic flows on manifolds of negative
sectional curvature.

There are a lot of works concerning the analysis of the  errors terms in (\ref{eq:1.1}) for different classes of dynamical systems as well as 
under different restrictions on the type of primitive closed trajectories considered (see \cite{kn:PP}, \cite{kn:PS2}, \cite{kn:PS3}, \cite{kn:PS4}, \cite{kn:L2}, \cite{kn:An} and the references there). For example, if $\varphi_t$ satisfies an approximative condition 
related to three primitive periods, Pollicott and Sharp \cite{kn:PS3} showed that there exists $\eta > 0$ such that
$$ \pi(x) = \frac{e^{h_T x}}{h_T x} \Bigl( 1 + \oo\Bigl(\frac{1}{x^{\eta}}\Bigr)\Bigr) \quad , \quad x \to +\infty\,.$$

On the other hand, for geodesic flows on negatively curved surfaces Pollicott and Sharp \cite{kn:PS2} proved a much sharper asymptotic: 
\begin{equation} \label{eq:1.2}
\pi(x) = {\rm li}\:(e^{h_T x}) + \oo(e^{cx}) \quad ,\quad 0 < c < h_T \; ,
\end{equation} 
where ${\rm li}\:(y) = \int_2^{y} \frac{1}{\log u} du.$ This results is based on estimates of the dynamical zeta function 
derived from strong spectral estimates for the iterations of the Ruelle transfer operator \cite{kn:D}. Recently it was shown 
that (\ref{eq:1.2}) holds for more general dynamical systems for which  strong spectral estimates for 
Ruelle transfer operators  were established (see \cite{kn:St2}, \cite{kn:St3}, \cite{kn:St5} and Sections 7- 9 below).

The purpose of this paper is to examine the number of closed trajectories with primitive periods lying in exponentially shrinking intervals
\be \label{eq:1.3}
(x - e^{-\delta x},\: x + e^{-\delta x})
\ee
as $x \to \infty$,  where $0 < \delta < h_T$. This question is closely related to the asymptotic behavior of sums of the form
\begin{equation} \label{eq:1.4}
\sum_{\sigma^n x = x} \psi_n (f^n(x)) \:\: ,\:\: n \to \infty\;,
\end{equation} 
where $\psi_n(t)$ ($t\in \R$) are functions with exponentially small support as $n \to \infty$, $f > 0$ is the so called 
roof function related to a given symbolic coding of the flow,  $\sigma$ is the shift in the corresponding symbol space, and 
$f^n(x) = f(x) + f(\sigma x) +...+ f(\sigma^{n-1} x)$. This type of ergodic sums for hyperbolic flows have been studied by many 
authors in the case when $\psi_n$ is the characteristic  function of an interval of the form ${\bf 1}_{[a\sqrt{n}, b\sqrt{n}]}$ 
(central limit theorems), ${\bf 1}_{[a, b]}$ or ${\bf 1}_{[z + p \epsilon_n, z + q \epsilon_n]}$ with $\epsilon_n \to 0$ not very fast 
(see \cite{kn:DP}, \cite{kn:L1}, \cite{kn:PS5}).  Moreover, in these works one assumes that $\int f d\nu = 0$, $\nu$ being a 
probability measure invariant with respect to $\sigma.$ In what follows below we simply replace $f$ by $g = f - \int f d\nu$.

In this paper we deal with functions of the form 
$\psi_n(t) = {\bf 1}_{[z + \alpha n + p \epsilon_n, z + \alpha n + q \epsilon_n]}(t)$ with $p < q,\:\epsilon_n = e^{-\delta n},\: \delta > 0$.  
To obtain an asymptotic for (\ref{eq:1.4}), we apply {\it strong spectral estimates} of the form (\ref{eq:1.5}) for the iterations of Ruelle 
transfer operators  (see \cite{kn:D}, \cite{kn:St2}, \cite{kn:St3} and Sects. 6-9).  On the other hand, the estimate 
(\ref{eq:1.2}) is based on the analytic continuation of the dynamical zeta function $Z(s)$ for 
$s_0 - \epsilon < \Re s \leq s_0,\: \epsilon >0$, $s_0$ being the abscissa of absolute convergence of $Z(s)$, and this continuation 
is obtained exploiting again the estimates (\ref{eq:1.5}). The second problem we deal with concerns the asymptotic of the number of 
primitive closed orbits. To obtain such an asymptotic is more difficult than estimating the number of periodic points of shifts maps 
in abstract symbol spaces. Clearly in this case  one has to estimate rather carefully the number of iterated periodic orbits involved in 
(\ref{eq:1.4}). 

Strong spectral estimates of the form (1.5) are known to hold for hyperbolic flows on basic  sets under certain additional regularity 
assumptions concerning the stable and unstable laminations over the basic set and under a local  nonintegrability condition 
(LNIC) ( see Sect. 9).  The latter appears to be a rather weak condition and is expected to be satisfied in most
(if not all) interesting cases. Indeed, it is already known that this condition is satisfied for contact Anosov flows (\cite{kn:St3}),
and for open billiard flows in $\R^n$ satisfying a certain additional regularity  assumption (\cite{kn:St4}).
In the present paper we show that (LNIC) always holds for arbitrary  basic sets of geodesic flows on hyperbolic manifols of 
constant negative curvature (see Lemma 3 in Sect. 7 which has an independent interest since it implies (\ref{eq:1.2})). 
In the latter case the stable/unstable laminations are smooth, so no extra regularity assumptions are necessary and the estimates (1.5) always hold.

To describe our results precisely we need to introduce some notation and definitions.
Let $\kappa_0 \geq 2$ be an integer and let $A = (A(i,j))_{i,j=1}^{\kappa_0}$
be a $\kappa_0\times \kappa_0$ aperiodic matrix of $0$'s and $1$'s.
Consider the  one-sided symbol space 
$$\saa = \{  (i_j)_{j=0}^\infty : 1\leq i_j \leq \kappa_0, A(i_j,\: i_{j+1}) = 1
\:\: \mbox{ \rm for all } \: j \geq 0\; \}$$
with the corresponding shift map $\sigma : \saa \longrightarrow \saa$, and 
for $\theta \in (0,1)$ let $\ff_\theta(\saa)$ be the space of
$d_\theta$-Lispchitz complex-valued functions on $\saa$ with the norms 
$\|\cdot \|_\infty$, $|\cdot |_\theta$ and $\|\cdot \|_\theta = \|\cdot \|_\infty + |\cdot |_\theta$
(see Sect. 2 for details). For a real-valued $g\in \ff_\theta(\saa)$ let $\Pr(g)$ be
the topological pressure of $g$ with respect to $\sigma$ (see Sect. 2).
Then there exists a unique $P_g \in \R$ such that $\Pr(-P_g\, g) = 0$. Consider the Ruelle operator
$$(\lc_{g} v)(\xi) = \sum_{\sigma \eta = \xi} e^{g(\eta)} v(\eta) \quad ,\quad \xi \in \ssa\;, \:v \in C(\ssa)\;.$$
When $g\in  \ff_\theta(\saa)$, this operator preserves the space $\ff_\theta(\saa)$ and it is bounded 
with respect to each of the norms $\|\cdot \|_\infty$ and $\|\cdot \|_\theta $. We will denote by 
$\|\lc_g \|_\infty$ and $\|\lc_g \|_\theta$ the norm 
of the operator $\lc_g$ with respect to any of these, respectively. Apart from that, given a real-valued function  $f\in \ff_\theta(\saa)$ and $a,u\in \R$ with $u \neq 0$, the operator $L_{(a + \i u)f}$ 
is bounded on $\ff_\theta(\saa)$ with respect to the norm
$$\|v\|_{\theta,u} = \|v\|_\infty + \frac{|v|_\theta}{|u|} \quad , \quad v\in \ff_\theta(\saa).$$
Throughout the paper we will need the following 

\ms

\begin{defn} We will say that the Ruelle transfer operators related to a
real-valued function $f\in \ff_\theta(\saa)$ 
are {\it weakly contracting} if for every $\ep > 0$ there exist constants $a_0 > 0$,
$\rho \in (0,1)$ and $A > 0$ $($possibly depending on $f$ and $\ep)$ such that
\be \label{eq:1.5}
\|\lc^n_{(-P_f+\i u)\,f}\|_{\theta, u} \leq   A \, \rho^n |u|^{\epsilon} \quad , \quad |u| \geq a_0\;,
\ee
for all integers $n \geq 0$.
\end{defn}
\ms

The above property is similar to the so called {\it strong spectral estimates for
Ruelle operators} related to basic sets of hyperbolic flows which we discuss in Sect. 9 below. 
There we also describe the conditions under which it is known that these estimates hold.

In the following we assume that $f(x) > 0$ for all $x \in \saa$ and set 
$$d_0 = \min_{x \in \saa}f(x),\: d_1 = \max_{x \in \saa} f(x).$$
Let $m_0$ be the equilibrium state of $- P_f \,f$. Then we have
$$\Pr(-P_f \, f) = h(m_0) - P_f\, \int f d m_0 = 0\, ,$$
where $h(m_0)$ is the measure-theoretic entropy of $m_0$ with respect to $\sigma$. 

Assuming that $f$ is non-lattice (see Section 2), there exists $\sigma_0 > 0$ such that
$$\frac{d^2 \Pr(- P_f f + \i\, u\,f )}{d u^2}\Big \vert_{u = 0} = - \sigma_0^2\;$$
(see \cite{kn:PP}). Set
$$\alpha = \alpha_f = \int f\, dm_0\;.$$

Given a constant $\delta > 0$, let
\be \label{eq:1.6}
\ep_n = e^{-\delta\, n} \quad , \quad n =1,2, \ldots\;.
\ee
For $p < q$ and an integer $n \geq 1$ set
$$I(z, p, q; \epsilon_n) = \# \{ \xi \in \saa:\: \mbox{\rm there exists}\:  
m \in \N\: \:\mbox{\rm with}\:\:\sigma^m(\xi) = \xi\:\: \mbox{\rm and} \:\: \:$$
$$ 
z + n\alpha + p \epsilon_n \leq f^m(\xi) \leq  z + n \alpha + q \epsilon_n\}\;.$$

Our first main result in this paper is the following

\begin{thm} Assume that the real-valued function  $f\in \ff_\theta(\saa)$ is non-lattice and
the Ruelle transfer operators related to $f$ are weakly contracting. 
Let $\epsilon_n = e^{-\delta n}$ with $0 < \delta < -\frac{\log \rho}{3}$, where 
$0 < \rho_1 < \rho < 1$ is such that $(\ref{eq:1.5})$ holds and $0 < \rho_1 < 1$ is the constant from Lemma $2$
in Sect. $3$ below. Then for any $0 \leq z \leq \alpha$ and any $p < q$ we have
\begin{equation} \label{eq:1.7}
\# \{ \xi \in \saa:\: \sigma^n(\xi) = \xi \:, \: z + n\alpha + p \epsilon_n \leq f^n(\xi) 
\leq z + n \alpha + q \epsilon_n\} \sim e^{P_f(z + n \alpha)} \frac{(q-p) \epsilon_n}{\sqrt{2 \pi} 
\sigma_0 \sqrt{n}}
\end{equation}
as $n \to \infty$, uniformly with respect to $z$.
\end{thm}
Here the notation $A(n) \sim B(n)$ as $n \to \infty$ means that $\lim_{n \to \infty} \frac{A(n)}{B(n)} = 1$
or equivalently $A(n) = B(n) (1 + o(1))$ with $o(1) \to 0$ as $n \to \infty.$ We also prove the following

\begin{thm}
Under the assumptions of Theorem $1$, assume that 
$\epsilon_n = e^{-\delta n}$ with $0 < \delta < - \frac{(\log \rho)\alpha}{3 d_1}.$ 
Then for any $0 \leq z \leq \alpha$, any $p < q$ and any fixed 
$a > 0$, setting $r = \frac{\pi}{4\alpha}$, we have
\begin{eqnarray} \label{eq:1.8}
e^{P_f(z + n \alpha)} (q-p) \epsilon_n \frac{1}{\sqrt{ \pi n} \sigma_0 }\frac{2 r}{a} \Bigl(1 + 
o_a(1)\Bigr) \leq I(z, p, q; \epsilon_n)\nonumber\\
\leq e^{P_f (z + n \alpha)} (q-p) \epsilon_n \frac{2 \sqrt{2 n}}{\sqrt{\pi} \sigma_0} \Bigl[ \sqrt{\frac{\alpha}{d_0}} - \sqrt{\frac{\alpha}{d_1}}  + o(1) \Bigr] 
, \quad n \to \infty\; .
\end{eqnarray}
uniformly with respect to $z$.
\end{thm}
The notation $o_a(1)$ means that the term $o_a(1)$ goes to 0 as $n \to \infty$ but the convergence to 0 depends on $a$.
As a simple consequence of the above results one obtains

\begin{cor} Under the assumptions of Theorem $2$ we have
\begin{equation} \label{eq:1.9}
\lim_{n \to \infty}   \frac{1}{n} \log I(z, p, q; \epsilon_n) = h(m_0) - \delta \;,
\end{equation}
therefore
$$\lim_{\delta \to 0} \lim_{n \to \infty}   \frac{1}{n} \log I(z, p, q; \epsilon_n) = h(m_0)\;.$$
\end{cor}

\begin{rem} 
If in the left-hand-side of $(\ref{eq:1.8})$ we take formally $\frac{r}{a} = \sqrt{2}n$, then we would have the 
same order with respect to $n$ in the right-hand-side and in the left-hand-side of $(1.8)$. However, this leads to 
a remainder $o_a(1)$ for which we have no control as $n \to \infty$ since $a$ depends of $n$. In this direction the 
result of Theorem $1$ is sharper, since we study the summation  only over the periodic points $x$ with $\sigma^n x = x$.
\end{rem}

The above results have natural consequences for hyperbolic flows. Here we state explicitly some
of them. Given a smooth flow $\varphi_t : M \longrightarrow M$ on a Riemannian manifold $M$, denote
by $\pp$ the set of all primitive  closed orbits of the flow, and for $\gamma\in \pp$ let 
$T_\gamma$ be the period (length) of $\gamma$.

Let $X = \H^{n+1}/\Gamma$ be a hyperbolic manifold, where
$\Gamma$ is a convex cocompact Kleinian group of transformation in $\H^{n+1}$,
and let $\varphi_t : M = S^*(X) \longrightarrow M$ be the geodesic flow on the unit 
cosphere bundle of $X$ (see Sect. $7$). Fix a Markov family $\rr = \{R_i\}_{i=1}^k$ for 
$\varphi_t$ over $\mt$ consisting of  rectangles $R_i = [U_i ,S_i ]$ such that the corresponding roof function
$\tau$ is non-lattice, set $U = \cup_{i = 1}^k {U_i}$ and let $\sigma : U \longrightarrow U$ be the naturally
defined shift map (see Sect. 6 for details).  Let $P = P_\tau \in \R$ be such that $\Pr(-P\, \tau) = 0$, where
$\Pr$ is the topological pressure with respect to $\sigma$, and let $m_0$ be the equilibrium state of 
$- P \,\tau$. Set $\di \alpha = \int_U \tau\, dm_0$, and let $\sigma_0 > 0$ be such that
$\frac{d^2 \Prf(- P \tau + \i\, u\,\tau )}{d u^2}\Big \vert_{u = 0} = - \sigma_0^2$. 
We then prove the following

\begin{thm} 
There exists $0 <\rho < 1$ such that  for every $0 < \delta < -\frac{\log \rho}{3}$, setting $\epsilon_n = e^{-\delta n}$
($n\in \N$), for any  $0 \leq z \leq \alpha$ and any $p < q$ we have
\begin{equation} \label{eq:1.10}
\# \{ x \in U :\: \sigma^n(x) = x \:, \: z + n\alpha + p \epsilon_n \leq \tau^n(x) 
\leq z + n \alpha + q \epsilon_n\} \sim e^{P(z + n \alpha)} \frac{(q-p) \epsilon_n}{\sqrt{2 \pi} 
\sigma_0 \sqrt{n}}
\end{equation}
as $n \to \infty$, uniformly with respect to $z$. Moreover, an analogue of $(\ref{eq:1.8})$ holds with $f$ replaced by $\tau$, 
$\Sigma_A^+$ by $U$ and $d_0 = \min_{x \in U} \tau(x),\: d_1 = \max_{x \in U} \tau(x).$
\end{thm}

The constant $\rho$ in the above theorem is such that the estimate (6.2) holds with $\rho$ and $1/\gamma < \rho < 1$ 
(see Sect. 6 for the definition of $\gamma$). Notice that for $x \in U$ with $\sigma^n(x) = x$, $\tau^n(x)$ is the length of a periodic trajectory passing through $x$, if $n$ is the smallest integer for which we have $\sigma^n (x) = x$. Thus we may derive a lower bound of the number of periodic trajectories with primitive periods lying in the interval
$(z + n\alpha + p \epsilon_n, z + n \alpha + q \epsilon_n)$ as $n \to \infty$ and we deduce a more precise result applying (\ref{eq:1.8}). 

A  similar result holds for general hyperbolic flows over basic sets satisfying some regularity conditions (see Sect. 6).

Our final result concerns open billiard flows in the exterior of several strictly convex domains
$K_1, \ldots, K_{\kappa_0}$ in $\R^N$, $N \geq 2$, satisfying some additional regularity conditions concerning 
the unstable and stable laminations through the non-wandering set (see Theorem 4 below and Sect. 8).
Since this flows has a natural coding by using boundary components,
in this case we get in a natural way results similar to Theorems 1 and 2 above. 
Namely, let $A$ be the $\kappa_0\times \kappa_0$ matrix such that $A(i,j) = 1$ if $i \neq j$ and
$A(i,j) = 0$ otherwise. Consider the space $\sa$ of double sequences with entries in
$\{ 1, \ldots, \kappa_0\}$ (see Sect. 2). Given any
$\xi = (\xi_j)_{j=-\infty}^\infty$ there exists a unique billiard trajectory $\gamma(\xi)$
in $\Omega = \overline{\R^N \setminus (K_1\cup \ldots \cup K_{\kappa_0})}$ with reflection
points $P_j(\xi) \in \partial K_{\xi_j}$. Set $f(\xi) = \|P_1(\xi)-P_0(\xi)\|$, and let
$m_0$ be the equilibrium state of $- P_f \,f$, where $\Pr(-P_f \, f)= 0$. Set
$\alpha = \alpha_f = \int f\, dm_0$. As before, there exists $\sigma_0 > 0$ such that
$$\frac{d^2 \Pr(- P_f f + \i\, u\,f )}{d u^2}\Big \vert_{u = 0} = - \sigma_0^2\;.$$
Finally,  let $\pp$ be the set of primitive closed billiard trajectories and let $\pp_n \subset \pp$ be the set of those 
primitive closed billiard trajectories $\gamma$ that have
exactly $n$ reflection points. Let $\ep_n$ be defined by (1.6). Set
$$I_n(z,p,q;\ep_n) = \# \{ \gamma \in \pp_n:\: z + n\alpha + p \en 
\leq T_{\gamma} \leq z + n \alpha + q \en\} \;.$$

Then we have the following

\begin{thm}  
Assume that the billiard flow $\varphi_t$ over its non-wandering set $\mt$ has regular 
distortion along unstable manifolds, satisfies the condition {\rm (LNIC)} and the local holonomy maps 
along stable laminations through $\mt$ are uniformly Lipschitz $($see Sect. $9$ below$)$. 
Then there exists $\delta_0$ such that
for $0 < \delta < \delta_0$ for any $0 \leq z \leq \alpha$, and any $p < q$ we have
\be \label{eq:1.11}
I_n(z,p,q;\ep_n) \sim e^{P_f( z + n \alpha)} 
\frac{ (q-p) \en }{\sqrt{2 \pi}n \sqrt{n} \sigma_0} \quad ,\quad n \to \infty\;,
\ee
uniformly with respect to $z$. Moreover,  for any fixed $a > 0$, setting $r = \frac{\pi}{4\alpha}$,
there exist  constants $C_0 > 0$ and $C_1 > 0$ such that
\begin{eqnarray} \label{eq:1.12}
e^{P_f( z + n \alpha)} (q-p) \en \frac{C_0 r}{a n\sqrt{ \pi n} \sigma_0}
\Bigl(1 + o_a(1)\Bigr) \nonumber \\
\leq  \#\{\gamma \in {\mathcal P}:\: z + n\alpha + p\epsilon_n \leq T_{\gamma} \leq z + n \alpha 
+ q \epsilon_n\} \nonumber\\
\leq  e^{P_f( z + n \alpha)} (q-p) \en \frac{1}{\sqrt{ \pi n} \sigma_0}\Bigl[C_1 \sqrt{\alpha} 
+ o(1)\Bigr] \quad , \quad n \to \infty \;.
\end{eqnarray}
In particular, the asymptotic $(\ref{eq:1.11})$ and the estimate $(\ref{eq:1.12})$ always hold when $\varphi_t$ 
satisfies the pinching condition {\rm(P)} over $\mt$  $($see Sect. $8)$.
\end{thm}

In fact, the condition (LNIC) (stated in Sect. 9 below)  follows from the result in \cite{kn:St4}
assuming that the local holonomy maps  along stable laminations through $\mt$ are $C^1$ (not just Lipschitz).
The latter is always the case if the pinching condition (P) (stated in Sect. 8 below) is satisfied.
As in Theorem 3, the constant $\delta_0 > 0$ in the above theorem depends on $\rho \in (0,1)$ and 
$1/\gamma \in (0,1)$ (see Sect. 6). 

A similar result holds for  other hyperbolic flows for which estimates similar to (\ref{eq:1.5}) are satisfied (see Sect. 9).

Sect. 2 contains  a few basic definitions from symbolic dynamics. Sects. 3, 4 and 5 are the main sections
in this paper -- they are devoted to the proofs of Theorems 1 and 2.  In Sect. 6 we consider general flows
over basic sets (satisfying certain additional conditions) and prove analogues of Theorems 1 and 2 -- see
Theorem 5 there. The proof of the latter is essentially a repetition of the arguments used in Sects. 3-5
with minor changes only.  Theorem 3 is derived as a consequence of Theorem 5 and Lemma 3 established in Sect. 7. Theorem 5  is also used in Sect. 8, where we prove Theorem 4.  In Sect. 9 we describe the main
result in \cite{kn:St3} concerning strong spectral estimates for Ruelle transfer operators which is used
essentially in Sects. 6-8.

\bs

\noindent
{\bf Acknowledgments.}
\footnotesize
Most of this work was accomplished during our stay at the Centre Interfacultaire Bernoulli, EPFL, Lausanne, 
as part of the Program `Spectral and Dynamical Properties of Quantum Hamiltonians'. Thanks are due to the
organizers of this Program and to the staff of the Centre Bernoulli for their hospitality and support.
We also grateful to Richard Sharp for the information he provided to us concerning Lemma 1 below as well as for his remarks on the previous version of the paper.
\normalsize

\section{ Preliminaries}
\renewcommand{\theequation}{\arabic{section}.\arabic{equation}}
\setcounter{equation}{0}

Let $\kappa_0 \geq 2$ be an integer and let $A = (A(i,j))_{i,j=1}^{\kappa_0}$
be a $\kappa_0\times \kappa_0$ matrix of $0$'s and $1$'s which is {\it aperiodic},
i.e. $A^M$ has strictly positive entries for some $M \geq 1$.
Consider the  symbolic space 
$$\sa = \{  (i_j)_{j=-\infty}^\infty : 1\leq i_j \leq \kappa_0, A(i_j, i_{j+1}) = 1
\:\: \mbox{ \rm for all } \: j\; \},$$
with the product topology and the {\it shift map} $\sigma : \sa \longrightarrow \sa$ 
given by $\sigma ( (i_j)) = ( (i'_j))$, where $i'_j = i_{j+1}$ for all $j$.
Given $0 < \theta < 1$, consider  the {\it metric} $d_\theta$ on $\sa$ defined by
$d_\theta(\xi,\eta) = 0$ if $\xi = \eta$ and $d_\theta(\xi,\eta) = \theta^m$ if
$\xi_i = \eta_i$ for $|i| < m$ and $m$ is maximal with this property.

In a similar way one deals with the one-sided subshift of finite type
$$\saa = \{  (i_j)_{j=0}^\infty : 1\leq i_j \leq \kappa_0, A(i_j, i_{j+1}) = 1
\:\: \mbox{ \rm for all } \: j \geq 0\; \},$$
where the {\it shift map} $\sigma : \saa \longrightarrow \saa$ is defined in
a similar way: $\sigma( (i_j)) = ( (i'_j))$, where $i'_j = i_{j+1}$
for all $j \geq 0$. The metric $d_\theta$ on $\saa$ is defined as above.
Let $\pi : \sa \longrightarrow \saa$ be the {\it natural projection}.

Let $B(\saa)$ be the space of bounded functions $g : \saa\longrightarrow \C$.
Given $f\in B(\saa)$ the {\it Ruelle transfer operator}
$\lc_f : B(\saa) \longrightarrow B(\saa)$ is defined by
$$\lc_fg(\xi) = \sum_{\sigma(\eta) = \xi} e^{f(\eta)}\, g(\eta)
\quad, \quad \xi \in \saa\;.$$

Let $\ff_\theta(\saa)$ denote the space of $d_\theta$-Lipschitz functions 
$g : \saa \longrightarrow \C$ 
with the norm
$\| f\|_\theta = \|f\|_\infty + |f|_\theta\;,$
where 
$$|f|_\theta = \sup \left\{ \frac{|f(\xi) - f(\eta)|}{d_\theta(\xi,\eta)} : 
\xi, \eta\in \saa\;, \; \xi \neq \eta\right\}\;.$$
If $f \in \ff_\theta(\saa)$, then  $\lc_f$ preserves the space $\ff_\theta(\saa)$.

\begin{defn} 
We say the function $f(x)$ on $\saa$ is non-lattice if there do not exist constants
$\gamma_0$ and $\gamma_1$, a function $G \in B(\saa)$ and an integer-valued function $Z \in B(\saa, \N)$ so that
$$f(x) = (G \circ \sigma)(x) - G(x)  + \gamma_0 + \gamma_1 Z(x),\: \forall x \in \ssa. $$
\end{defn}

Given a real-valued $F \in \ff_\theta(\saa)$ the {\it topological pressure} $\Pr\,(F)$ of 
$F$ is defined by
$$\Pr (F) = \sup_{\mu \in {\mathcal M}_{\sigma}} 
\left[h(\mu) + \int_{\saa} F\, d\mu \right]\,,$$
where ${\mathcal M_{\sigma}}$ is the set of all probability measures on $\saa$ invariant with 
respect to $\sigma$ and $h(\mu)$ is the {\it measure-theoretic entropy} of $\mu$ 
with respect to $\sigma$ (see e.g. \cite{kn:PP} for more details). Notice that for small $|u|$ we may define the
pressure $\Pr(F + \i u G)$ for real-valued functions $F, G \in \ff_{\theta}(\ssa)$ since the Ruelle operator ${\mathcal L}_{F}$ has a simple "maximal" eigenvalue (see Section 4 and Proposition 4.7 in \cite{kn:PP}).

\section{The case of a Markov shift }
\renewcommand{\theequation}{\arabic{section}.\arabic{equation}}
\setcounter{equation}{0}

\subsection{Representations of $S(n)$}

Let  $\sigma: \saa \longrightarrow \saa$ be the shift on $\saa$ 
and let the real-valued function $f \in \ff_\theta(\saa)$ for some $\theta \in (0,1)$.

Assume that $f$ is non-lattice and that the Ruelle transfer operators related to 
$f$ are weakly contracting, so that (\ref{eq:1.5}) holds.
As in Sect. 1, let $P = P_f$ be such that $\Pr(-P\, f) = 0$, and let
$m_0$ be the equilibrium state of $- P \,f$ so that
$$ h(m_0) - P\, \int f \, d m_0 = \Pr(-P \, f) = 0\, .$$
Below we will write simply $P$ instead of $P_f$ since the function $f$ is fixed in our considerations.

Set $\alpha = \int f \, dm_0$ and consider a sequence 
$\{\epsilon_n\}_{n \in \N},\:\epsilon_n > 0,\: \epsilon_n \to 0$ such that 
\begin{equation} \label{eq:3.1}
\epsilon_n  = e^{- \delta n},
\end{equation}
with $0 < \delta < -\frac{\log \rho}{3}$, where $\rho \in ( \rho_1 , 1)$ is the constant that appears in 
(\ref{eq:1.5}) and $0 < \rho_1 < 1$ is 
the constant from Lemma 2 below.
Let $\chi: \R \longrightarrow \R^+$ be a $C^k$ ($k \geq 3$), function with compact support. 
Set $\chi_n(x) = \chi(\epsilon_n^{-1}(x- z))$ and
$g = f - \int f \, d m_0$, and note that $\int g \, d m_0 = 0.$

We will study the behavior of
$$S(n) := \sum_{x\in \saa\, ,\,\sigma^n x = x} \chi_n(g^n(x)) $$
$$= \frac{1}{2 \pi} \int_{-\infty}^{\infty} \Bigl(\sum_{\sigma^n x = x} e^{\i u g^n(x)}\Bigr) \hat{\chi}_n(u) du,$$
where $\hat{\chi}_n(u) = e^{-\i z u}\epsilon_n \hat{\chi}(\epsilon_n u)$ and $\hat{\chi}(u)$ is the Fourier transform of $\chi.$  Introduce the function
$\omega_n(y) = e^{-\xi y}\chi_n(y)$ with $\xi = -P$. 
Then
$$\hat{\omega}_n(u) = \int e^{-\i u y} e^{-\xi y} \chi_n(y) dy = \hat{\chi}_n(u - \i \xi) 
= \epsilon_n e^{-\i zu} e^{-\xi z}\hat{\chi}(\epsilon_n(u - \i\xi)),$$
and  
\begin{eqnarray*}
S(n) 
& = & \ssn e^{\xi g^n(x)}\omega_n(g^n(x)) 
= \frac{1}{2 \pi} \int_{-\infty}^{\infty} \Bigl ( \ssn e^{(\xi + \i u) g^n(x)} \Bigr) 
\hat{\omega}_n(u) du\\
& = & \frac{ e^{P(z + n \alpha)} \epsilon_n}{2 \pi}\int_{-\infty}^{\infty} \Bigl(\ssn e^{-P f^n(x) + \i ug^n(x)}\Bigr) 
e^{-\i zu}\hat{\chi}(\epsilon_n (u - \i \xi)) du\\
& = & \frac{e^{ P (z + n\alpha)} \epsilon_n}{2 \pi}\Bigl[\int_{|u| < a} + \int_{a \leq |u| \leq c} + \int_{|u| > c}\Bigr],
\end{eqnarray*}
where $a > 0$ will be chosen sufficiently small and $c \gg 1$ sufficiently large. With this partition we have 
$S(n) = I_{1,n} + I_{2, n} + I_{3, n}.$

For periodic points we have the following Lemma which follows from the fact that $\Pr(-P f) = 0$ 
and the proof of the statement (ii) of  Theorem 5.5 in \cite{kn:PP}.

\begin{lem} 
There exists $0 < \theta_1 < 1$ and $a > 0$ such that for $|u| \leq a$  we have
$$\ssn e^{-P f^n(x) + \i ug^n(x)} = e^{n \Prf(-P f + \i u g)} + \oo(n \theta_1^n).$$
\end{lem}

Our choice of $\xi = -P$ implies that $\xi g = - P\, f+ P\, \alpha$ and
$${\mathcal L}^n_{(\xi + \i u)g} = e^{P\, n \alpha} {\mathcal L}^n_{- P\, f + \i u g}\;.$$ 

Next we have
$$\frac{d \Pr(-P\, f + \i u g)}{du}\Big\vert_{u = 0} = \i \int g \,d m_0 = 0\;.$$
Moreover, since  $g$ is non-lattice, we deduce
$$\frac{d^2 \Pr(- P\, f + \i ug)}{d u^2}\Big \vert_{u = 0} = - \sigma_0^2$$
for some $\sigma_0 > 0.$

The representation of the sum 
$$\ssn e^{-P f^n(x) + \i ug^n(x)} = e^{-\i n \alpha}\ssn e^{(- P + \i u) f^n(x)}$$
for $|u| > a$ is more complicated and we will use the so called Ruelle's Lemma in the form
proved in \cite{kn:W}. Let $\chi_i$ be the characteristic function of the cylinder
$${\mathcal C}_i =\{ \eta \in \Sigma_A^+: \eta_0 = i\} \quad ,\quad i = 1,...,\kappa_0\;.$$
Fix an arbitrary point $x_i \in {\mathcal C}_i$. Then we have the following

\begin{lem} 
There exists a constant $\rho_1 \in (0,1)$ such for $a_0 > 0$ and $b_0 > 0$ 
and every $\epsilon > 0$ there exist constant $C_{\epsilon} > 0$ so that for $|t| \leq a_0, \:|u| \geq b_0$ we have the estimate
\begin{equation} \label{eq:3.2}
\Bigl|\ssn e^{(t + \i u) f^n(x)} - \sum_{i = 1}^{\kappa_0} {\mathcal L}_{(t + \i u) f}^n \chi_i(x_i)\Bigr|
\leq C_{\epsilon} |u| \sum_{m = 2}^n 
\Bigl(\|{\mathcal L}_{(t + \i u )f}^{n-m}\|_{\theta} \, \rho_1^m \, e^{m (\epsilon + \Prf(t f))}\Bigr)
\end{equation}
for all $n \in \N$.
\end{lem}

This lemma  was proved in \cite{kn:W} generalizing and completing some points of the proof of a similar lemma 
in \cite{kn:PS2} and \cite{kn:N} proved for surfaces and $C^1$ regular foliations. In our case we treat manifolds with arbitrary dimensions and (\ref{eq:3.2}) is established in \cite{kn:W} for functions $f \in {\mathcal F}_{\theta}(\Sigma_A^+).$
 Notice that in the setting of Sect. 6 we can choose $\rho_1 = 1/\gamma$, where $\gamma$ is as in (6.1).

\section{ Asymptotic of $S(n)$}
\renewcommand{\theequation}{\arabic{section}.\arabic{equation}}
\setcounter{equation}{0}

\subsection{Asymptotic for $|u| < a$}

We start with the analysis of $I_{1,n}.$
Choosing $a > 0$ sufficiently small and changing the coordinates on $(-a, a)$ to $v = \frac{\sigma_0 u}{\sqrt{2}}$, 
we can write
$$e^{\Prf(- P f + \i u g)} = (1 - v^2 + \i Q(v))\;,$$
where $Q(v)$ is a real-valued function such that $Q(v) = \oo(|v|^3)$
(see Lemma 1.2 (3) in \cite{kn:PS5} and Proposition 2.2 in \cite{kn:PS6}).

Modulo terms involving $\oo(n\theta_1^n)$, the term $I_{1, n}$ has the form
$$I_{1, n} =  e^{P (z + n\alpha)}\frac{\epsilon_n\sqrt{2}}{2 \pi\sigma_0}\int_{-b}^b e^{-\i u(v) z}
\Bigl[(1 - v^2 + \i Q(v))^n  \Bigr]\hat{\chi}(\epsilon_n(u(v) + \i P))
dv\;$$
with $b = \frac{\sigma_0 a}{\sqrt{2}}$ and $u(v) = \frac{\sqrt{2} v}{\sigma_0}.$
We have $\hat{\chi}(\epsilon_n(u(v) + \i \,P)) = \hat{\chi}(0) + \epsilon_n \oo( 1 +|v|)$ 
and $e^{- \i u(v) z} = 1 + \oo_z (|v|).$ The leading term of $I_{1, n}$ becomes
$$e^{P (z + n \alpha)}\frac{\epsilon_n \sqrt{2}\hat{\chi}(0)}{2 \pi\sigma_0} 
\int_{-b}^b(1 - v^2)^n dv 
= e^{P (z + n \alpha)}\frac{\epsilon_n \sqrt{2}\hat{\chi}(0)}{2 \pi\sigma_0}
\int_0^{b^2} \frac{(1 - w)^n}{w^{1/2}} dw\;.$$

Next
$$\int_0^{b^2} \frac{(1 - w)^n}{w^{1/2}} dw = \int_0^1 \frac{(1-w)^n}{w^{1/2}}dw + 
\oo((1 - b^2)^n) \sim \frac{\sqrt{\pi}}{\sqrt{n}}\;.$$
as $n \to +\infty$. Here we use the formula 
$$\int_0^1 (1 - w)^n w^{q/2 -1} dw = \frac{\Gamma(n+1)\Gamma(q/2)}{\Gamma(n+1 + \frac{1}{2}q)}$$
and apply the Stirling approximation for $\Gamma(m).$

On the other hand, 
$$|(1 - v^2 + \oo(|v|^3))^n - (1 - v^2)^n| \leq 
\mbox{\rm Const}\, \sum_{j = 1}^n C^n_j (1- v^2)^{n-j} a^j |v|^{3j}$$
and we can estimate the integral of the right-hand-side by $\oo(\frac{1}{n})$. We refer 
to \cite{kn:PS1} for the details of this calculation. The integration of the perturbation $\oo(n\theta_1^n)$ yields a negligible term and we conclude that
\begin{equation} \label{eq:4.1}
I_{1, n} \sim e^{P (z + n \alpha)}\frac{\epsilon_n \hat{\chi}(0)}{\sqrt{2 \pi}\sigma_0 \sqrt{n}}
\quad ,\quad n \to +\infty.
\end{equation}
Notice that $\hat{\chi}(0) = \int  \chi(y) dy > 0.$\\

\subsection{Asymptotic for $a \leq |u| \leq c$}

First consider the integral
$$J_{2, n} =  e^{P (z + n \alpha)}\frac{\epsilon_n}{2 \pi} \int_{a < |u| 
\leq c}e^{- \i u z}\sum_{i = 1}^{\kappa_0}\lc_{- P f + \i u g}^n \chi_i(x_i) du$$
with $c \gg 1$ sufficiently large which will be chosen below. 

Notice that
$$\lc_{-P f + \i u g}^n = e^{-\i n u \alpha} \lc_{-P f + \i u f}^n.$$
Since $-P f$ is non-lattice,  for $0 < a \leq |u| \leq c$ the operator 
${\mathcal L}_{-P f + \i u f}$ has no eigenvalues $\mu, \:|\mu| = 1$ (see for instance \cite{kn:PP}) 
and the spectral radius of ${\mathcal L}_{-P f + \i u f}$ is strictly less than 1. Thus, there exist
$\beta = \beta(a, c),\: 0 < \beta < 1$ and $C_{a, c} > 0$ such that we have
\begin{equation} \label{eq:4.2}
\|{\mathcal L}^n_{-P f + \i u f}\|_{\theta} \leq C_{a, c}\beta^n\quad , \quad a \leq |u|\leq c,\: \forall n \in \N\;.
\end{equation}
On the other hand,
\begin{equation} \label{eq:4.3}
|\hat{\chi}(\epsilon_n(u - \i \xi))| \leq C_m \frac{e^{c_0 |\epsilon_n \xi|}}{\epsilon_n^m |u|^m}
\quad ,\quad |u| \geq a,\: m \in \N,\:m \leq k\;,
\end{equation}
with $c_0 > 0$ depending on the support of $\chi.$ 
Using (\ref{eq:4.2}) and (\ref{eq:4.3}) with 
$k = 0$, for large $n$ we get 
$$|J_{2, n}| \leq C_{a, c,\chi} e^{P (z + n \alpha)} e^{c_0 P}\frac{\epsilon_n }{ 2 \pi}\beta^n 
\int_{a \leq |u| \leq b} du \leq C'(a, c, \chi)  e^{P (z + n \alpha)} \epsilon_n \beta^n.$$
Next to estimate the sum in the right hand side of (\ref{eq:3.2}) we choose $\epsilon $ small and we increase
$0 < \beta < 1$ , if necessary so that $\frac {\rho_1 e^{\epsilon}}{\beta} =  \theta_2 < 1.$ Therefore,

$$\sum_{j=2}^n \|{\mathcal L}_{-Pf + \i u f}^{n-j}\|_{\theta} (\rho_1 e^{\epsilon})^j \leq C_{a, c}\beta^n\sum_{j = 2}^n \Bigl(\frac{\rho_1 e^{\epsilon}}{\beta}\Bigr)^j \leq C'_{a, c} \beta^n$$
and we repeat the argument for the estimation of $J_{2,n}.$
Finally, we get
\begin{equation} \label{eq:4.4}
I_{2, n} = {\mathcal O} \Bigl( e^{P (z + n \alpha)}\frac{\epsilon_n}{n}\Bigr).
\end{equation}

\subsection{Asymptotic for $|u| > c$}
We apply lemma 2 with $t = -P$. In this case $\Pr(-Pf) = 0$ and we must examine
$$I_{3, n} = e^{P (z + n \alpha)}\frac{\epsilon_n}{2 \pi} \int_{|u| > c} 
e^{- \i u (z + n \alpha)}\hat{\chi}(\epsilon_n(u - \i \xi))\Bigl[\sum_{i = 1}^{\kappa_0}\lc_{(- P  + i u)f}^n \chi_i(x_i)$$
$$ + \oo_{\epsilon}\Bigl(|u| \sum_{j = 2}^n \|\lc_{(-P + \i u) f}^{n-j}\|_{\theta} (\rho_1 e^{\epsilon})^j\Bigr)\Bigr] du = J_{n, 3} + R_{n, 3}\;.$$
It follows from (\ref{eq:1.5}) that if $c$ is large enough we have for $|u| > c$ and every $\nu > 0$ the estimates

\begin{equation} \label{eq:4.5}
\|\lc_{(-P + \i u) f}^n\|_{\infty} + \frac{|\lc_{(-P + \i u)f}^n |_{\theta}}{|u|} \leq A_{\nu} \rho^n |u|^{\nu}, 
\: \forall n \in \N.
\end{equation}
We choose $\eta > 0$ and $\nu > 0$ small enough in order to arrange $-\delta \geq \frac{\log \rho}{3 + \nu} + \eta.$ Then 
$$\epsilon_n^{3 + \nu} = e^{-\delta (3 + \nu) n} 
\geq \rho^n e^{\eta(3 + \nu) n}.$$
For the sum over $i = 1,...,\kappa_0$ we apply (\ref{eq:4.5}) with $\nu$ to estimate the $\|.\|_{\infty}$ norm and  for large $n$ we get
$$A_{\nu}\rho^n \int_{|u| > c} |u|^{\nu} |\hat{\chi}(\epsilon_n (u - i\xi))| du 
\leq A_{\nu} e^{-\eta(3 + \nu) n}\int_{|u| > c}  \epsilon_n^2 |u|\hat{\chi}(\epsilon_n (u - i\xi))|du$$
$$\leq \frac{1}{n} \int |y| |\hat{\chi}(y - i \epsilon_n \xi)| dy = \oo\Bigl(\frac{1}{n}\Bigr)\;.$$
The integral involving $\oo_{\epsilon}\Bigl(|u| \sum_{m = j}^n \|\lc_{(-P + \i u) f}^{n-j}\|_{\theta} (\rho_1 e^{\epsilon})^j\Bigr)$
is dealt in the same way. First for $\epsilon$ small we arrange the inequality
$$\frac{\rho_1 e^{\epsilon}}{\rho} = \theta_3 < 1$$
increasing, if it is necessary, $0 < \rho < 1$ in (\ref{eq:4.5}).
Then we have
$$C_{\epsilon} |u| \sum_{j= 2}^n \|\lc_{(-P + \i u) f}^{n-j}\|_{\theta} (\rho_1 e^{\epsilon})^j \leq 
C_{\epsilon, \nu} |u|^{2 + \nu} \rho^n \sum_{j= 2}^n \theta_3^j \leq C_{\epsilon, \nu}' |u|^{2 + \nu} \rho^n.$$
Consequently,
$$|R_{n, 3}| \leq e^{P(z + n \alpha)}\frac{\epsilon_n}{2 \pi}C_{\epsilon, \nu}'\int_{|u| \geq c} \rho^n |u|^{2 + \nu}|\hat{\chi}(\epsilon_n (u - i\xi))|du$$
$$\leq e^{P(z + n \alpha)}\frac{\epsilon_n}{2 \pi}C_{\epsilon, \nu}'e^{-\eta(3 + \nu)n}\int_{|u| \geq c} \epsilon_n (\epsilon_n |u|)^{2+\nu} |\hat{\chi}(\epsilon_n (u - i\xi))|du$$
and for large $n$ we get
$$C_{\epsilon, \nu}' e^{-\eta(3 + \nu)n} \int |y|^{2+ \nu} |\hat{\chi}(y - \i \epsilon_n \xi)| dy = \oo\Bigl(\frac{1}{n}\Bigr).$$
Thus, we conclude that 
$$I_{3, n} = {\mathcal O} \Bigl( e^{P (z + n \alpha)}\frac{\epsilon_n}{n}\Bigr)\,.$$
Consequently, for $n \to +\infty$ we obtain the following

\begin{prop}  Let $f$ be non-lattice and such that the Ruelle transfer operators related to
$f$ are weakly contracting. Let $\epsilon_n = e^{-\delta n}$, where $0 < \delta < -\frac{\log \rho}{3}$ with $0 < \rho_1 < \rho < 1$ such that $(\ref{eq:4.5})$ holds.  Then
\begin{equation} \label{eq:4.6}
S(n) \sim  e^{P (z + n \alpha)}\frac{\epsilon_n }{\sqrt{2 n \pi}\sigma_0} 
\int \chi(y) dy,\: n \to +\infty.
\end{equation}
\end{prop}

Now it is easy to pass from $\chi \in C_0^k(\R)$ to an indicator function ${\bf 1}_{[p, q]}$ of the 
interval $[p, q]$ repeating the argument in \cite{kn:PS5}. For completeness we give the proof. Given 
$\eta > 0$, choose cut-off functions $\chi^{-},\: \chi^{+}$ so that 
$\chi^{-} \leq {\bf 1}_{[p, q]} \leq \chi^{+}$ and
$$q - p - \eta \leq \int \chi^{-}(x) dx \leq \int \chi^{+} (x) dx \leq q - p + \eta\, .$$
Using (4.6), we get
$$ \limsup _{n \to +\infty} \frac{\sqrt{n}}{e^{n h(m_0)} \epsilon_n} 
\#\{x \in {\rm Fix}(\sigma^n):\:  z + n\alpha + p \epsilon_n \leq f^n(x) \leq z + n \alpha + q \epsilon_n\}$$
$$\leq \limsup _{n \to +\infty} \frac{\sqrt{n}}{e^{n h(m_0)} \epsilon_n} 
\ssn \chi_n^{+}(g^n(x)) \leq e^{P z}\frac{q - p + \eta}{\sqrt{2 \pi} \sigma_0},$$
$$\liminf_{n \to +\infty} \frac{\sqrt{n}}{e^{n h(m_0)} \epsilon_n}
\#\{x \in {\rm Fix}(\sigma^n):\:  z + n\alpha + p \epsilon_n \leq f^n(x) \leq z + n \alpha + q \epsilon_n\}\;$$
$$\geq \liminf_{n \to +\infty} \frac{\sqrt{n}}{e^{n h(m_0)} \epsilon_n} 
\ssn \chi_n^{-}(g^n(x))\geq e^{P z}\frac{q - p - \eta}{\sqrt{2 \pi} \sigma_0}\;.$$
Since $\eta > 0$ is arbitrary, we deduce that for any $z \in \R$ we have
\be \label{eq:4.7}
\#\{ x \in {\rm Fix}\: (\sigma^n):\: z + n \alpha + p \epsilon_n \leq f^n(x) \leq z + n\alpha + q \epsilon_n\} \sim e^{P(z + n\alpha)}\frac{(q - p)\epsilon_n}{\sqrt{2 \pi}\sigma_0 \sqrt{n}}.
\ee
Moreover, the asymptotic is uniform for $z$ in a compact interval. This proves 
Theorem 1. \endofproof

\bs

To study the distribution of primitive periods we need to examine the function
$$S_{min}(n) := \sum_{\sigma^n x = x,\: n \:{\text minimal }}\: 
{\bf 1}_{[z + p\epsilon_n, z + q \epsilon_n]}(g^n(x))\;, $$
where the summation is over all points $x \in \ssa$ such that 
$n = \min\{ m \in \N:\:\sigma^m x = x\}$.
For this purpose observe that
$$\#\{ x \in {\rm Fix}\: (\sigma^n):\: z + n \alpha + p \epsilon_n 
\leq f^n(x) \leq z + n\alpha + q \epsilon_n\}$$
$$ = S_{min}(n) + \sum_{\sigma^m x = x,\:m\: {\rm minimal}\atop n/m = k \in \N,\: k \geq 2} 
{\bf 1}_{[z + p\epsilon_n, z + q\epsilon_n]}(g^{km}(x))= S_{min}(n) + S_r(n)\;.$$

Any $x \in {\rm Fix}\: (\sigma^n)$ defines a {\it periodic $n$-orbit}
$\gamma = \{ \sigma^j(x) : 0 \leq j \leq n-1\}$. We will say that $\gamma$ is 
{\it primitive} if $n\geq 1$ is the smallest integer with $\sigma^n(x) = x$.
The number $T_\gamma = f^n(x)$ will be called the {\it $f$-period} of $\gamma$.
Let $\pp_n$ be the {\it set of all primitive periodic $n$-orbits}. 

For the $f$-periods of primitive periodic $n$-orbits $\gamma \in {\mathcal P}_n$,
we must divide by $n$ since $\gamma$ contains $n$ points in ${\rm Fix}\:(\sigma^n)$. 
Thus, by (4.7) we get an upper bound
\begin{equation} \label{eq:4.8}
\#\{ \gamma \in {\mathcal P}_n:\: z + n \alpha + p \epsilon_n \leq T_{\gamma} 
\leq z + n\alpha + q \epsilon_n\} 
\leq e^{P(z + n\alpha)}\frac{(q - p)\epsilon_n}{\sqrt{2 \pi}\sigma_0 n\sqrt{n}}(1 + o(1)), 
n \to \infty.
\end{equation}
To obtain a lower bound, we assume that $h(m_0) > \delta$ where $h(m_0) = P\alpha$.
Notice that
$$S_{min}(n) \geq e^{P (z + n \alpha)}\frac{\epsilon_n }{\sqrt{2 n \pi}\sigma_0}
\Bigl( (q-p)  - o(1)\Bigr) - S_r(n),\: n \to +\infty.$$
Thus it is sufficient to have an upper bound for $S_r(n).$ Consider a term in $S_r(n)$ having the form
$$G_{n, m} = \sum_{\sigma^m x = x,\:m\: {\rm minimal},\atop n/m = k, k \geq 2}
{\bf 1}_{[z + p\epsilon_n, z + q \epsilon_n]}(g^{km}(x)).$$
with some fixed divisor $m \in \N$ of $n$.
Then for $\sigma^m x = x$ we get $g^{km}(x) = kf^m(x) - km\alpha$ and for large $n$ we have
$$ \frac{z}{k} +  m \alpha - \epsilon_m \leq \frac{z}{k} + m \alpha +\frac{p}{k} \epsilon_n     
\leq    f^m(x) 
\leq \frac{z}{k} + m \alpha + \frac{q}{k} \epsilon_n \leq \frac{z}{k} + m \alpha + \epsilon_m\;.$$
We choose $\eta > 0$ so that $0 < \eta <  h(m_0) - \delta.$
Next we fix an integer $k_0 \in \N$ such that 
\begin{equation} \label{eq:4.8}
0 < \frac{h(m_0)}{k_0} < h(m_0) - \delta - \eta.
\end{equation}
Consider two cases: $(i)\: 2 \leq k \leq k_0,\: (ii)\: k > k_0.$ In the case $(i)$ we have 
$m = \frac{n}{k} \geq \frac{n}{k_0} \longrightarrow \infty$ as $n \to \infty.$
Since
$${\bf 1}_{[z + p \epsilon_n, z + q \epsilon_n]}(g^{km}(x)) \leq {\bf 1}_{[\frac{z}{k} 
- \epsilon_m, \frac{z}{k} + \epsilon_m]}(g^m(x))\;,$$
we can apply (\ref{eq:4.6}) with $m$ replaced by $n$ and $p = -1,\: q = 1$. Thus,
$$
G_{n,m} \leq e^{P(z/k + m\alpha)}\frac{2\epsilon_m}{\sqrt{2 \pi}\sigma_0 \sqrt{m}}(1 + o(1))
\leq e^{P z} e^{(h(m_0) - \delta)\frac{n}{2}}\frac{\sqrt{2 k_0}}{\sqrt{\pi}\sigma_0 \sqrt{n}}
(1 + o(1))$$
$$\leq e^{P z}e^{(h(m_0) - \delta)n} e^{-\mu n} \frac { \sqrt{ 2k_0}}{\sqrt{ \pi} \sigma_0 \sqrt{n}} 
( 1 + o(1)),$$
where $0 < 2\mu < h(m_0) - \delta.$
Summing over $2 \leq k \leq k_0$, we obtain
\begin{eqnarray} \label{eq:4.9}
\sum_{2 \leq k \leq k_0} G_{n, m} \leq e^{P(z + n \alpha)} \epsilon_n e^{-\mu n}\frac{(k_0 -1)
\sqrt{ 2k_0} }{ \sqrt{\pi} \sigma_0 \sqrt{n}} (1 +o(1))\\
\leq B_{\sigma_0, k_0} e^{P n \alpha} \epsilon_n \frac{e^{-\mu n}}{\sqrt{n}}(1 + o(1)).\nonumber
\end{eqnarray}

Passing to the case $(ii)$, notice that for $k > k_0$ we cannot guarantee that $m = n/k$ goes 
to $+ \infty$ as $n \to \infty$. For this reason we apply a crude estimate for the  number of  
$f$-periods  $T_{\gamma}$ of periodic rays $\gamma$. Namely, since $f$ is non-lattice, we exploit
the estimate for the number of primitive periods 
$$\#\{T_{\gamma} \leq x\} \sim \frac{e^{h_T x}}{h_T x},\: x \to \infty,$$
where $h_T > 0$ is the topological entropy of the suspended symbolic flow related to $f$ 
(see \cite{kn:PP}). This estimate is based on the analysis of the behavior of the following {\it dynamical zeta function}
$$Z(s) = \sum_{n=1}^{\infty} \frac{1}{n}\sum_{\sigma^n x = x}  e^{-s f^n(x)}$$
(see for instance, \cite{kn:PP}). Notice that in our case the abscissa of absolute convergence of $Z(s)$ is exactly
$h_T = P \alpha.$  Thus, for $0 \leq \frac{z}{k_0} \leq \frac{\alpha}{k_0},\: q \epsilon_n \leq q$ we get
$$\# \{ \gamma:\: T_{\gamma} \leq \frac{z}{k} + \frac{n}{k}\alpha + q \epsilon_n\} \leq C_{q, k_0}
\frac{e^{P\frac{n}{k_0}\alpha}}{P \frac{n}{k_0} \alpha}(1 + o(1))\quad ,\quad n \to \infty\;.$$

Summing over $k_0 < k \leq n/2$ and taking into account (\ref{eq:4.8}), we obtain 
\begin{equation} \label{eq:4.10}
\sum_{k_0 < k \leq n/2} G_{n,m} \leq k_0 C_{z,q, k_0}\frac{e^{h(m_0) \frac{n}{k_0}}} {2 h(m_0)} 
(1 + o(1))\leq A_{z, q, k_0} e^{P n \alpha} \epsilon_n e^{-n\eta} (1 + o(1)),\: n \to \infty.
\end{equation}
Combining (\ref{eq:4.9}) and (\ref{eq:4.10}), we deduce
$$S_{min}(n) \geq e^{P (z + n \alpha)} \frac{\epsilon_n}{\sqrt{2 \pi n} \sigma_0} ((q-p) - o(1)).$$
Finally, in $S_{min}(n)$ every periodic primitive orbit is counted $n$ times and we obtain the asymptotic (\ref{eq:1.11}) for $0 < \delta < \min\{-\frac{\log \rho}{3}, h(m_0)\}.$

\begin{rem} The analysis in this section follows the approach in \cite{kn:PS5}, Section $4$. However, 
the argument in \cite{kn:PS5} exploits Lemma $4.2$ there which is not proved and it seems that in that form 
the lemma is not correct. Our arguments are based on Lemmas $1$ and $2$ above. Moreover, the 
investigation of the case $|u| > c$ with exponentially decreasing $\epsilon_n$ seems impossible without using
strong spectral estimates of the form $(\ref{eq:1.5})$.

\end{rem}

\section{Asymptotic of $S(n, m)$}
\renewcommand{\theequation}{\arabic{section}.\arabic{equation}}
\setcounter{equation}{0}

In this section we study the counting function of primitive periodic orbits related to $\sigma^m x = x$ and 
having $f$-periods in the interval $[z + n\alpha +p \epsilon_n,\: z +  n\alpha + q\epsilon_n], \: 0 \leq z \leq \alpha.$  Let
$$d_0 = \min_{x \in \ssa} f(x),\:\: d_1 = \max_{x \in \ssa} f(x).$$
The non-lattice condition on $f$ implies $d_0 < d_1.$ We assume in this section that $\delta$ in (\ref{eq:3.1}) satisfies
\begin{equation} \label{eq:5.1}
0 < \delta < -\frac{(\log \rho) \alpha}{3 d_1}
\end{equation}

If $\gamma$ is a primitive periodic orbit with $m$ points such that 
$z + n\alpha + p \epsilon_n \leq T_\gamma \leq z + n\alpha + q \epsilon_n$, then
$m d_0 \leq T_{\gamma} \leq m d_1$ and 
\begin{equation} \label{eq:5.2}
\frac{n\alpha}{d_1} + \oo(\epsilon_n) \leq  \frac{ z + n\alpha +p \epsilon_n}{d_1} 
\leq  m \leq \frac{ z + n\alpha + q\epsilon_n}{d_0}  \leq \frac{(n+1)\alpha}{d_0} + \oo (\epsilon_n).
\end{equation}
Introduce the function
$$\chi_{n, m}(x) = \chi\Bigl((x -z - n\alpha + m \alpha)\epsilon_n^{-1}\Bigl)$$
and note that 
$\hat{\chi}_{n,m}(u) = e^{- \i z u}e^{-\i (n - m)\alpha u} \epsilon_n \hat{\chi}(\epsilon_n u).$ 
Next consider the sum
$$S(n, m) = \sum_{\sigma^m x = x} \chi_{n, m}(g^m(x)).$$

Using the notation of the previous section, we get
$$S(n, m) = \frac{\epsilon_n}{2 \pi} \int_{-\infty}^{\infty} 
\Bigl( \sum_{\sigma^m x = x} e^{(\xi + \i u)g^m(x)}\Bigr)e^{-\i (n - m)\alpha u} 
e^{- \i u z - \xi z}e^{-\xi(n-m)\alpha} \hat{\chi}(\epsilon_n(u - \i \xi)) du.$$

We consider three zones of integration: $|u| \leq a, \: a < |u| \leq c,\: |u| \geq c$, where 
$a$ is small enough and $c$ is sufficiently large.

Repeating the argument of Section 4, we must study for sufficiently small $b > 0$ the integral
$$I_{n, m} = \int_{-b}^b \Bigl[ (1 - v^2 + \i Q(v))^m 
e^{- \i \frac{\sqrt{2}}{\sigma_0} v} e^{-\i (n-m)\alpha  
\frac{\sqrt{2}}{\sigma_0}v} \hat{\chi}(\epsilon_n(\frac{\sqrt{2}}{\sigma_0}v - \i\xi)) dv.$$
The only difference is the presence of the oscillatory factor
\begin{equation} \label{eq:5.4}
e^{-\i (n-m)\alpha \frac{\sqrt{2}}{\sigma_0}v} .
\end{equation}

Notice that the leading term becomes
$$\hat{\chi}(0)\int_{-b}^{b} (1 - v^2)^m e^{-i(n -m)\alpha \frac{\sqrt{2}}{\sigma} v} dv 
= 2\hat{\chi}(0)\int_0^b (1 - v^2)^m \cos\Bigl((n-m) \alpha \frac{\sqrt{2}}{\sigma_0}v \Bigr) dv$$
for $m$ satisfying (\ref{eq:5.2}). We will obtain a {\it lower bound} for number of the 
periods taking into account only the $f$-periods of periodic orbits related to $\sigma^m x = x$ for which
$$|n - m| \leq \frac{\pi }{4 \alpha  a}.$$
For such $m$ we get
$$\int_0^b (1 - v^2)^m \cos\Bigl((n-m) \alpha \frac{\sqrt{2}}{\sigma_0}v \Bigr)dv 
\geq \frac{1}{\sqrt{2}} \int_0^b (1 - v^2)^m dv $$
$$ = \frac{1}{\sqrt{2}}\Bigl[ \int_0^1(1- v^2)^m dv - \int_b^1(1 - v^2)^m dv\Bigr] $$
$$ \sim \frac{\sqrt{\pi}}{\sqrt{2 m}} + \oo\Bigl((1 - b^2)^m\Bigr),\: m \to \infty.$$
Setting 
$r = \frac{\pi}{4 \alpha },$
we have
$$\sqrt{\frac{\pi}{2}}\sum_{n - \frac{r}{a} \leq m \leq n + \frac{r}{a}} \frac{1}{\sqrt{m}} 
\sim \sqrt{2 \pi}\Bigl[\Bigl(n + \frac{r}{a}\Bigr)^{1/2} - \Bigl( n - \frac{r}{a}\Bigr)^{1/2} \Bigr]$$
$$= \sqrt{2 \pi n} \Bigl[ \Bigl(1 + \frac{r}{an}\Bigr)^{1/2} - \Bigl( 1 - \frac{r}{an}\Bigr)^{1/2}\Bigr] = \frac{\sqrt{2 \pi }r}{a \sqrt{n}} + {\mathcal O}_a\Bigl(\frac{1}{n \sqrt{n}}\Bigr),\: n \to \infty.$$

On the other hand,
$$\sum_{n - \frac{r}{a} \leq m \leq n + \frac{r}{a}} (1- b^2)^m \leq C(r, a) \exp(n\log (1- b^2)).$$
To obtain an upper bound for $I_{m, n}$, note that
$$\sum_{\frac{n\alpha}{d_1}-1 \leq m \leq \frac{(n + 1)\alpha}{d_0}+1} \frac{\sqrt{\pi}}{\sqrt{m}} 
\leq \sqrt{\pi}\int_{\frac{n \alpha}{d_1}-2} ^{\frac{(n+1)\alpha}{d_0} + 2} x^{-1/2}dx$$
$$\leq 2 \sqrt{\pi n} \Bigl( \sqrt{\frac{\alpha}{d_0}} -\sqrt {\frac{\alpha}{d_1}} 
+ \oo\Bigl(\frac{1}{n}\Bigr)\Bigr) ,\: n \to \infty.$$
The analysis of lower order terms goes without any change and we obtain 
\begin{eqnarray} \label{eq:5.5}
\frac{e^{P(z+ n \alpha)} \epsilon_n \hat{\chi}(0)}{\sqrt{\pi n} \sigma_0} \frac{2 r}{a} 
\Bigl( 1 + {\mathcal O}_a \Bigl(\frac{1}{n \sqrt{n}}\Bigr)\Bigr) \leq  I_{n,m} \\ \nonumber
\leq \frac{e^{P (z+ n \alpha)} \epsilon_n \hat{\chi}(0) 2 \sqrt{2 n}}{\sqrt{\pi} \sigma_0} 
\Bigl( \sqrt{\frac{\alpha}{d_0}} -\sqrt {\frac{\alpha}{d_1}}  + {\mathcal O}\Bigl( \frac{1}{n}\Bigr)\Bigr).
\end{eqnarray}

The integral over $b \leq |v| \leq c$ can be treated as in Section 4 since we have a factor
$\epsilon_n$ and $m \to \infty.$ In the analysis of the integral over $|v| > c$ we must take 
into account that for the operators $\lc_{-Pf + \i u f}^{m-j}$ in (\ref{eq:3.2}) the estimate
(\ref{eq:1.5}) gives a decay with $\rho^{m-j}$ and not with $\rho^{n-j}.$ On the other hand, 
$m \geq \frac{n\alpha}{d_1} + \oo(\epsilon_n)$ and the assumption (\ref{eq:5.1}) imply
$\rho^{m} \leq C \rho^{\frac{n\alpha}{d_1}} \leq \epsilon_n^3 e^{-\eta n}$
with some $\eta > 0.$ Thus the analysis in the previous section goes without change and the 
integral over $|v| > c$ yields negligible terms.

To pass to an indicator function, we exploit the same argument as in the previous section to get
\begin{eqnarray} \label{eq:5.5}
e^{P(z + n \alpha)} \epsilon_n \frac{(q-p)}{\sqrt {\pi n}\sigma_0} \frac{2r}{a} \Bigl(1 + o_a(1)\Bigr) \leq I(z, p, q; \epsilon_n)\\
\leq e^{P ( z +n \alpha)} \epsilon_n \frac{(q -p) 2 \sqrt{2 n}}{\sqrt{\pi} \sigma_0} 
\Bigl( \sqrt{\frac{\alpha}{d_0}} -\sqrt {\frac{\alpha}{d_1}}  + o(1)\Bigr),\: n \to \infty. \nonumber
\end{eqnarray}
This completes the proof of Theorem 2.\endofproof

\bs

Now we pass to the analysis of the counting function
$$S_{min}(n, m): = \sum_{\sigma^m x = x, m\: {\rm minimal}} {\bf 1}_{[z + n\alpha 
+ p \epsilon_n, z + n\alpha + q \epsilon_n]}(f^m(x)). $$
As in the previous section we write $S_{min}(n, m) = S(n, m) - S_r(n, m)$ and we will find an upper
bound of $S_r(n, m).$ To do this, we will apply an argument similar to that used in Section 4 and 
we sketch below the necessary modifications. Let 
$$\sigma^s(x) = x,\:s \: {\rm minimal},\:  m = ks,\: s \in \N, k \geq 2\;,$$
and let $n = kt + l,\: t \in \N,\: 0 \leq l \leq k-1.$ Then $f^m(x) = kf^{s}(x)$ and  
$z + n\alpha + p \epsilon_n \leq f^m(x) \leq z + n\alpha + q \epsilon_n$ for large $n$ implies
$$\frac{z + l \alpha}{k} + t \alpha - \epsilon_t \leq f^s(x) \leq \frac{z + l \alpha}{k} 
+ t \alpha + \epsilon_t.$$
We consider two cases: $(i)\: 2 \leq k \leq k_0, \: (ii)\: k > k_0$ and we choose $k_0$ large 
enough in order to have
(\ref{eq:4.8}). In the case $(i)$, for fixed $m$ satisfying (\ref{eq:5.2}), we consider the divisors $s$ of $m$ with $\sigma^s(x) = x$ and we apply (\ref{eq:5.6}) with $p = - 1, q = 1$, replacing 
$n$ by $ t$ and $z$ by $\frac{z + l \alpha}{k}.$ This is possible since $t = \frac{n}{k} - \frac{l}{k} \geq \frac{n}{k_0} - 1 \to \infty$ as $n \to \infty.$ Thus we obtain an upper bound with order
$\oo\Bigl(e^{h(m_0)t} \epsilon_t\Bigl) \leq \oo\Bigl(e^{(h(m_0) - \delta)\frac{n}{2}}\Bigr)$ 
since $t \leq n/2.$ Then we repeat the argument in Section 4 and get a negligible term. For 
$k > k_0$ we apply again a crude estimate 
$$\#\{ \gamma:\: T_{\gamma} \leq \frac{z}{k} + \frac{n}{k} \alpha + \frac{q}{k}\epsilon_n\} \leq C_{z,q, k_0} \frac{e^{P \frac{n}{k_0} \alpha}}{P \frac{n}{k_0}\alpha} (1 + o(1))$$
and we exploit (\ref{eq:4.8}). Summing with over $2 \leq k \leq \frac{(n + 1) \alpha}{2d_0}$, 
and then over $m \leq \frac{(n +1)\alpha}{d_0}$ we obtain an upper bound for $S_r(n,m)$ and we
conclude that $S_r(n,m)$ yields a negligible term. Consequently, for $S_{min}(n,m)$ we deduce the
same estimates as in (\ref{eq:5.6}). For the counting function of the primitive periodic rays in 
$\pp$ we must divide the upper bound of $S_{min}(n,m)$ by $\frac{n\alpha}{d_1}$ and the lower bound of $S_{min}(n,m)$ by $\frac{(n + 1)\alpha}{d_1}.$ Thus we obtain the estimates (\ref{eq:1.12}).

\def\Intu{\Int^{u}}
\def\Ints{\Int^{s}}

\section{General hyperbolic flows over basic sets}
\renewcommand{\theequation}{\arabic{section}.\arabic{equation}}
\setcounter{equation}{0}

Let $\varphi_t : M \longrightarrow M$ be a $C^2$ Axiom A flow on a
$C^2$ complete (not necessarily compact) Riemannian manifold $M$. 
A $\varphi_t$-invariant closed subset $\mt$ of $M$ is called {\it hyperbolic} if $\mt$ contains
no fixed points  and there exist  constants $C > 0$ and $0 < \lambda < 1$ such that 
there exists a $d\varphi_t$-invariant decomposition 
$T_xM = E^0(x) \oplus E^u(x) \oplus E^s(x)$ of $T_xM$ ($x \in \mt$) into a direct sum of non-zero 
linear subspaces,
where $E^0(x)$ is the one-dimensional subspace determined by the direction of the flow
at $x$, $\| d\varphi_t(u)\| \leq C\, \lambda^t\, \|u\|$ for all  $u\in E^s(x)$ and $t\geq 0$, and
$\| d\varphi_t(u)\| \leq C\, \lambda^{-t}\, \|u\|$ for all $u\in E^u(x)$ and  $t\leq 0$. Here
$\|\cdot \|$ is the {\it norm} on $T_xM$ determined by the Riemannian metric on $M$.

A non-empty compact $\varphi_t$-invariant hyperbolic subset $\mt$ of $M$ which is not a single 
closed orbit is called a {\it basic set} for $\varphi_t$ if $\varphi_t$ is transitive on $\mt$ 
and $\mt$ is locally maximal, i.e. there exists an open neighbourhood $V$ of
$\mt$ in $M$ such that $\mt = \cap_{t\in \R} \varphi_t(V)$.  When $M$ is compact and $M$ itself
is a basic set, $\phi_t$ is called an {\it Anosov flow}.

Let $\mt$ be a basic set for $\varphi_t$.
For $x\in \mt$ and $\ep > 0$ sufficiently small, let 
$$W^s_\ep(x) = \{ y\in M : d (\varphi_t(x),\varphi_t(y)) \leq \epsilon \: \mbox{\rm for all }
\: t \geq 0 \; , \: d (\varphi_t(x),\varphi_t(y)) \to_{t\to \infty} 0\: \}\; ,$$
$$W^u_\ep(x) = \{ y\in M : d (\varphi_t(x),\varphi_t(y)) \leq \epsilon \: \mbox{\rm for all }
\: t \leq 0 \; , \: d (\varphi_t(x),\varphi_t(y)) \to_{t\to -\infty} 0\: \}$$
be the (strong) {\it stable} and {\it unstable manifolds} of size $\epsilon$. Then
$E^u(x) = T_x W^u_\ep(x)$ and $E^s(x) = T_x W^s_\ep(x)$.

Throughout this section we will assume that $\Lambda$ is a basic set for $\varphi_t$
such that the local holonomy maps along stable laminations through $\Lambda$ are
uniformly Lipschitz (see Sect. 9 below). Following \cite{kn:R1}, a subset $R$ of $\mt$ will be called a {\it rectangle} if it has the form
$$R = [U,S] = \{ [x,y] : x\in U, y\in S\}\;,$$ 
where $U$ and $S$ are admissible subsets of $W^u_\ep(z) \cap \mt$ and $W^s_\ep(z) \cap \mt$, 
respectively, for some $z\in \mt$ (cf. e.g. \cite{kn:D} or Sect. 2 in \cite{kn:St3}).
For such $R$, given $\xi = [x,y] \in R$, we will denote $W^u_R(\xi) = \{ [x',y] : x'\in U\}$ and
$W^s_R(\xi) = \{[x,y'] : y'\in S\} \subset W^s_{\ep'}(x)$. Denote by $\Intu(U)$ (resp.
$\Ints(S_i)$) the {\it interior} of the set $U$ in $W^u_\ep(z) \cap \mt$ (resp. $W^s_\ep(z) \cap \mt$)
and set $\Int(R) = [\Intu(U), \Ints(S)]$. Similarly, for $\xi = [x,y] \in R$ set $\Intu(W^u_R(\xi)) = [\Intu(U), y]$
and $\Ints(W^us_R(\xi)) = [x,\Ints(S)]$.

Let $\rr = \{ R_i\}_{i=1}^k$ be a family of rectangles with $R_i = [U_i  , S_i ]$,
$U_i \subset W^u_\ep(z_i) \cap \mt$ and $S_i \subset W^s_\ep(z_i)\cap \mt$, respectively, 
for some $z_i\in \mt$.   Set  $R =  \cup_{i=1}^k R_i\; .$
The family $\rr$ is called {\it complete} if  there exists $T > 0$ such that for every $x \in \mt$,
$\varphi_{t}(x) \in R$ for some  $t \in (0,T]$.  The {\it Poincar\'e map} $\pp: R \longrightarrow R$
related to a complete family $\rr$ is defined by $\pp(x) = \varphi_{\tau(x)}(x) \in R$, where
$\tau(x) > 0$ is the smallest positive time with $\varphi_{\tau(x)}(x) \in R$.
The function $\tau$  is called the {\it first return time} associated with $\rr$. Notice that  $\tau$ is
constant on each of the set $W_{R_i}^s(x)$, $x \in R_i$. 
A complete family $\rr = \{ R_i\}_{i=1}^k$ of rectangles in $\mt$ is called a 
{\it Markov family} of size $\chi > 0$ for the  flow $\varphi_t$ if $\diam(R_i) < \chi$ for all $i$ and: 

\ms

(a)  for any $i\neq j$ and any $x\in \Int(R_i) \cap \pp^{-1}(\Int(R_j))$ we have   
$$\pp(\Ints (W_{R_i}^s(x)) ) \subset \Ints (W_{R_j}^s(\pp(x)))\quad, \quad 
\pp(\Intu(W_{R_i}^u(x))) \supset \Intu(W_{R_j}^u(\pp(x)))\; ;$$

(b) for any $i\neq j$ at least one of the sets $R_i \cap \varphi_{[0,\chi]}(R_j)$ and
$R_j \cap \varphi_{[0,\chi]}(R_i)$ is empty.

\ms

The existence of a Markov family $\rr$ of an arbitrarily small size $\chi > 0$ for $\varphi_t$
follows from the construction of Bowen \cite{kn:B} (cf. also  Ratner \cite{kn:Ra}). 

Let $\rr = \{R_i\}_{i=1}^k$ be a Markov family for $\varphi_t$ over $\mt$.
Setting $U = \cup_{i=1}^k U_i$, the {\it shift map} $\sigma : U \longrightarrow U$
is defined by $\sigma = p\circ \pp$, where $p : R \longrightarrow U$
is the projection along the leaves of local stable manifolds. Let $\hR$ be the set of all
$x \in R$  whose orbits do not have common points with the boundary of $R$. Set $\hU = U\cap \hR$.
It is well-known (\cite{kn:B}) that $\hR$ is a residual subset of $R$ that has full measure
with respect to any Gibbs measure on $R$. The same applies to $\hU$ in $U$.

Denote by $C(U)$ the space of bounded continuous functions $h : U \longrightarrow \C$
with the usual norm $\|h\|_\infty = \sup_{x\in U} |h(x)|$. Notice that
$\tau$ is continuous on $\hU$, however in general $\tau$ could be discontinuous on $U$. 
Next, denote by $\clip(U)$ the {\it space of Lipschitz functions} $v : U \longrightarrow \C.$ 
For such $v$  let $\Lip(v)$ denote the {\it Lipschitz constant} of $v$, and for $u\in \R$, $u \neq 0$, define
$$\|v\|_{\lip,u} = \|v\|_\infty + \frac{\Lip(v)}{|u|} \quad , \quad
\|v\|_{\lip} = \|v\|_\infty + \Lip(v)\;.$$

\begin{rem}
The function $\tau$ is {\it locally Lipschitz} on $R$ in the following sense:
there exists a constant $\Lip(\tau) > 0$ such that if $x,y \in R_i$
for some $i$ and $\sigma(x), \sigma(y) \in R_j$ for some $j$, then
$|\tau(x) - \tau(y)| \leq \Lip(\tau)\, d(x,y)$. The map $\pp$ has a similar
property. Moreover, it is easy to see that for any $h\in \clip(U)$ and any $s\in \C$ 
the operator $L_{h+s\, \tau}$  preserves the space $\clip(U)$.
\end{rem}

The hyperbolicity of the flow on $\mt$ implies the existence of
constants $c_0 \in (0,1]$ and $\gamma_1 > \gamma > 1$ such that
\begin{equation}
c_0 \gamma^m\; d (u_1,u_2) \leq 
d (\sigma^m(u_1), \sigma^m(u_2)) \leq \frac{\gamma_1^m}{c_0} d (u_1,u_2)
\end{equation}
whenever $\sigma^j(u_1)$ and $\sigma^j(u_2)$ belong to the same  $U_{i_j}$ 
for all $j = 0,1 \ldots,m$.  

From now on we will assume that $\mt$ is a fixed basic set for $\varphi_t$
and $\rr = \{R_i\}_{i=1}^k$ is a fixed Markov family for $\varphi_t$ over $\mt$ consisting of 
rectangles $R_i = [U_i ,S_i ]$.
Let $\aa = (\aa_{ij})_{i,j=1}^k$ be the matrix given by $\aa_{ij} = 1$ 
if $\pp(\Int (R_i)) \cap \Int (R_j) \neq  \e$ and $\aa_{ij} = 0$ otherwise. 
It is well-known (\cite{kn:BR}) that the Markov family $\rr$ can be chosen so that
{\bf $\tau$ is non-lattice}. From now on we will assume that $\rr$ is chosen in this way.

Given a Markov family $\rr$, one defines a natural
symbol space $\Sigma_\aa$ and a natural map $\W : \sA \longrightarrow R$ such that
$\W\circ \sigma = \pp\circ \W$, where $\sigma : \Sigma_\aa \longrightarrow \Sigma_\aa$
is the shift map. However, in general $\W$ is not one-to-one and this presents
certain difficulties in trying to apply Theorems 1 and 2 to count numbers of periodic
orbits in $\mt$. Instead of using the symbol space $\sa$ and the coding map $\W$, 
here we just use the arguments from the proofs of Theorems 1 and 2 in a slightly
different setting to derive similar results.

Using the setup in Sect. 1, let $P = P_\tau \in \R$ be such that $\Pr(-P\, \tau) = 0$, where
$\Pr$ is the topological pressure with respect to $\sigma : U \longrightarrow U$,
and let $m_0$ be the equilibrium state of $- P \,\tau$. 
Since $\tau$ is non-lattice, there exists $\sigma_0 > 0$ such that
$\di \frac{d^2 P(- P \tau + \i\, u\,\tau )}{d u^2}\Big \vert_{u = 0} = - \sigma_0^2\;.$
As in Sect. 1, set $\di \alpha = \int_U \tau\, dm_0\;.$

In the present setting the analogue of Definition 1 reads the following.

\ms

\begin{defn} We will say that the Ruelle transfer operators related to a
real-valued function $f\in \clip(U)$ are {\it weakly contracting} if for every $\ep > 0$ 
there exist constants $a_0 > 0$, $\rho \in (0,1)$ and $A > 0$ (possibly depending on $f$ 
and $\ep$) such that
\be 
\|L_{(- P_f  + \i u) f}^n 1\|_{\lip, u} \leq A \, \rho^n |u|^{\epsilon} \quad , \quad |u| \geq a_0\;,
\ee
for all integers $n \geq 0$.
\end{defn}
\ms

Set $d_0 = \inf_{x\in U} \tau(x)$ and $d_1 = \sup_{x\in U} \tau(x)$.
The following theorem comprises the analogues of Theorems 1 and 2 
in the present setting.


\noindent
\begin{thm}
Assume that the Ruelle transfer operators related to $\tau$ are weakly contracting.

\ms

$(a)$ Let $\epsilon_n = e^{-\delta n}$ with $0 < \delta < -\frac{\log \rho}{3}$, where $0 < 1/\gamma < \rho < 1$ 
and {\rm (6.2)} holds with $\rho$. Then for any $0 \leq z \leq \alpha$ and any $p < q$ we have
\begin{equation}
\sharp \{ x \in U :\: \sigma^n(x) = x \:, \: z + n\alpha + p \epsilon_n \leq \tau^n(x) 
\leq z + n \alpha + q \epsilon_n\} 
\sim e^{P (z + n \alpha)} \frac{(q-p) \epsilon_n}{\sqrt{2 \pi} \sigma_0  \sqrt{n}}
\end{equation}
as $n \to \infty$, uniformly with respect to $z$.

\ms

$(b)$ Let $\epsilon_n = e^{-\delta n}$ with $0 < \delta < -\frac{(\log \rho)\alpha}{3d_1}$. 
Then for any $0 \leq z \leq \alpha$, any $p < q$ and any fixed 
$a > 0$, setting $r = \frac{\pi}{4\alpha}$, we have
\begin{eqnarray}
e^{P(z + n \alpha)} (q-p) \epsilon_n \frac{1}{\sqrt{ \pi n} \sigma_0 }\frac{2 r}{a} \Bigl(1 + 
o_a(1)\Bigr)\Bigr) \leq I(z, p, q; \epsilon_n)\nonumber\\
\leq e^{P (z + n \alpha)} (q-p) \epsilon_n \frac{2 \sqrt{2 n}}{\sqrt{\pi} \sigma_0} 
\Bigl[ \sqrt{\frac{\alpha}{d_0}} - \sqrt{\frac{\alpha}{d_1}}  + o(1) \Bigr] \quad , \quad n \to \infty\; ,
\end{eqnarray}
uniformly with respect to $z$.

\end{thm}
\noindent
{\it Proof.} This is a repetition of the arguments in Sects. 3, 4 and 5.
Here we give a very brief  sketch of these for completeness. In the present setting 
$R$ plays the role of $\sa$ and $U$ that of $\saa$. Moreover,  $f = \tau$,
which is constant on stable leaves of rectangles $R_i$ (i.e. $f$ depends on future coordinates only). In general, $\tau$
is not continuous on $U$, however as mentioned in Remark 2 above, $L_{s\tau}$ preserves the space
$\clip(U)$ for any $s\in \C$.

Next, set $\di g = \tau - \int_U \tau\, dm_0$, choose the function $\chi$ as in Sect. 3,
and define $\chi_n$, $S_n$ and $\omega_n$ as in Sect. 3.1, where 
$h_T = P$ is the topological entropy of $\varphi_t$ on $\mt$ and $\xi =  - h_T$.
Lemma 2 applies in the present setting without change, so we have the
inequality (3.2), as well with $\|\lc_{(-P_f + \i u)f}\|_{\theta}$ replaced by $\|\lc_{(-P_f + \i u)f}\|_{\lip}$. 
The argument at the end of Sect. 3 applies without change.

Next, the analytic arguments in Sect. 4 also apply without change 
and for $f = \tau$ the argument in Sect. 5 works without any change.
\endofproof

\section{Geodesic flows on manifolds of constant negative curvature}
\renewcommand{\theequation}{\arabic{section}.\arabic{equation}}
\setcounter{equation}{0}

Let $X$ be a complete (not necessarily compact) connected Riemannian manifold of constant curvature 
$K = -1$ and dimension $\dim(X) = n+1$, $n \geq 1$,  and let $\varphi_t : M = S^*(X) \longrightarrow M$ 
be the geodesic flow on the {\it unit cosphere bundle} of $X$.  According to a classical result of 
Killing and Hopf,
any such $X$ is a {\it hyperbolic manifold}, i.e. $X$ is isometric to $\H^{n+1}/\Gamma$, where 
$$\H^{n+1} = \{ (x_1, \ldots, x_{n+1}) \in \R^{n+1} : x_{1} > 0\}$$
is the upper half-space in $\R^{n+1}$ with the {\it Poincar\'e metric} 
$ds^2(x) = \frac{1}{x^2_{1}} (dx_1^2 + \ldots + dx_{n+1}^2)$
and $\Gamma$ is a {\it Kleinian group}, i.e. a discrete group of isometries (M\"obius transformations) 
of $\H^{n+1}$. See e.g. \cite{kn:Ratc}  for basic information on hyperbolic manifolds.  
Given a hyperbolic manifold $X = \H^{n+1}/\Gamma$,
the {\it limit set} $L(\Gamma)$ is defined as the set of accumulation points of all $\Gamma$ orbits in 
$\overline{\partial \H^{n+1}}$, the topological closure of $\partial \H^{n+1} = \{ 0\} \times \R^n$ 
including $\infty$.

Throughout this section we will assume that $\Gamma$ is torsion-free and finitely generated (then
$\Gamma$ is {\it geometrically finite}) and {\it non-elementary}, i.e. $L(\Gamma)$ is infinite
(then $L(\Gamma)$ is a closed non-empty nowhere dense subset of 
$\partial \overline{\H^{n+1}}$ without isolated points; see e.g. Sect. 12.1 in \cite{kn:Ratc}). 
A geometrically finite Kleinian group with no parabolic elements is called 
{\it convex cocompact}. If $X$ is compact, then $\Gamma$ is called a {\it cocompact lattice}.

The {\it non-wandering set}  $\mt$ of  $\varphi_t : M  \longrightarrow M$ 
(also known as the {\it convex core} of $X = \H^{n+1}/\Gamma$)
is the image in $M$ of the set of all points of $S^*(\H^{n+1})$ generating geodesics with end 
points in $L(\Gamma)$. When $\Gamma$ is convex cocompact, the non-wandering set $\mt$  is compact. 

From now on we will assume that $\Gamma$ is a non-elementary convex cocompact Kleinian group.

As in Sect. 6, let $\rr = \{R_i\}_{i=1}^k$ be a fixed Markov family for $\varphi_t$ over $\mt$ 
consisting of  rectangles $R_i = [U_i ,S_i ]$ such that the corresponding roof function
$\tau$ is non-lattice. Let $P = P_\tau \in \R$ be such that $\Pr(-P\, \tau) = 0$, where
$\Pr$ is the topological pressure with respect to $\sigma : U \longrightarrow U$,
and let $m_0$ be the equilibrium state of $- P \,\tau$. Set
$\di \alpha = \int_U \tau\, dm_0$, and let $\sigma_0 > 0$ be such that
$\frac{d^2 \Prf(- P \tau + \i\, u\,\tau )}{d u^2}\Big \vert_{u = 0} = - \sigma_0^2$.

\begin{lem}
The Ruelle transfer operators related to $\tau$ are weakly contracting.
\end{lem}

Now Theorem 3 follows from Lemma 3 and the arguments in Sect. 3-5 as we have obtained Theorem 5 in the previous section.

\bs

\def\tq{\tilde{q}}
\def\hq{\hat{q}}
\def\hLa{\hat{\Lambda}}
\def\tx{\tilde{x}}
\def\ty{\tilde{y}}
\def\tz{\tilde{z}}
\def\pH{\partial \H^{n+1}}
\def\txo{\tilde{x}^{(0)}}
\def\tso{\tilde{\sigma}_0}
\def\tPsi{\tilde{\Psi}}
\def\tmt{\tilde{\Lambda}}
\def\sh{S^*(\H^{n+1})}
\def\zoo{z^{(1)}}
\def\zo{z^{(0)}}
\def\xo{x^{(0)}}
\def\yoo{y^{(1)}}
\def\xoo{x^{(1)}}
\def\hy{\hat{y}}
\def\hx{\hat{x}}
\def\hz{\hat{z}}
\def\tg{\tilde{g}}
\def\tS{\tilde{S}}
\def\tsi{\tilde{\sigma}}
\def\ttt{\tilde{t}}
\def\la{\langle}
\def\ra{\rangle}

\noindent
{\it Proof of Lemma} 3.
We will use an argument from \cite{kn:St5}.
Let $p: \H^{n+1} \longrightarrow X  = \H^{n+1}/\Gamma$ and 
$\hat{p} : S^*(\H^{n+1}) \longrightarrow M = S^*(X)$ be the natural
projections. Consider the {\it geodesic flow} 
$\phi_t: S^*(\H^{n+1}) \longrightarrow S^*(\H^{n+1})$  on $\H^{n+1}$.
Recall that the geodesics in $\H^{n+1}$ are either straight lines perpendicular to 
$\partial \H^{n+1} = \{ x\in \R^{n+1} : x_{1} = 0\}$
or semi-circles with centers in $\partial \H^{n+1}$ whose planes are perpendicular to 
$\partial \H^{n+1}$.

It is known that  the {\it non-wandering set} $\Lambda\subset M$ of $\varphi_t$ has the form 
$\mt = \hat{p}(\hLa)$, where 
$\hLa$ is the set of those $x\in S^*(\H^{n+1})$ such that both $\lim_{t \to\infty} \phi_t(x)$ 
and $\lim_{t \to-\infty} \phi_t(x)$ belong to the {\it limit set} $L(\Gamma)$ of the group $\Gamma$.
The assumptions made above imply that $L(\Gamma)$ is a non-empty $\Gamma$-invariant
closed subset of $\partial \H^{n+1}$ without isolated points (see Ch. 12 in \cite{kn:Ratc}).

A {\it horosphere} in $\H^{n+1}$ is
either an $n$-sphere in $\overline{\H^{n+1}}$ tangent to $\partial \H^{n+1}$, or an $n$-plane in
$\H^{n+1}$ parallel to $\partial \H^{n+1}$. Let $S$ be a horosphere and $x\in S\cap \H^{n+1}$.
If $S$ is an $n$-sphere, denote by $\nu_S(x)$ the {\it outward normal} to $S$ at $x$
with $\| \nu_S(x)\| = 1/x_{1}$, while if $S$ is
an $n$-plane, set $\nu_S(x) = -\frac{1}{x_{1}}\, e_{1} = \frac{1}{x_{1}}\, (-1,0, \ldots, 0)$. 
The stable and unstable manifolds for $z = (x,\nu_S(x))$ in $S^*(\H^{n+1})$ are given by
$$W^s(x) = \{ (y, -\nu_S(y)) : y\in S\cap \H^{n+1}\}\quad, \quad 
W^u(x) = \{ (y,\nu_S(y)) : y\in S\cap \H^{n+1}\}\;,$$
so obviously the local stable and unstable foliations are smooth. The projections of the latter  via 
$\hat{p}$ give the local stable and unstable foliations in $M$.

It is also straightforward to check that $\varphi_t$ has uniform distortion along
unstable manifolds over $\mt$. 

To check this it is again enough to work on the universal cover $\H^{n+1}$.
Let $\tz = (z,\zeta) \in S^*(\H^{n+1})$ and $t > 0$.  Since the isometry group of $\H^{n+1}$ is both
point and direction transitive, we may assume that $z = (1,0, \ldots, 0)$ and $\zeta = -e_{1}$. 
Then 
$$W^u(\tz) = \{ (y,-e_{1}) : y_{1} = 1\} \quad, \quad
W^u(\phi_t(\tz)) = \{ (w,-e^{-t}\, e_{1}) :  w_{1} = 1- e^{-t}\}\;.$$
Obviously,  for any smooth curve $\gamma$ in $W^u(\tz)$ of length $\ell_\gamma$, the length of
$\phi_t(\gamma)$ is exactly $e^t \cdot \ell_\gamma$.  Thus, for any
$\tx,\ty  \in W^u(\tz)\setminus \{\tz\}$ we have 
$\frac{d (\tz,\tx)}{d (\tz, \ty)} =  
\frac{d (\phi_t(\tz),\phi_t(\tx))}{d (\phi_t(\tz), \phi_t(\ty))}$. 
Since $\hat{p}$ is a local isometry conjugating the geodesic flows $\varphi_t$ and $\phi_t$, 
it follows that $\varphi_t$ has uniform distortion along
unstable manifolds over $\mt$.

It remains to check the condition (LNIC) of Sect. 9 below.
Again we will work on the universal cover $\H^{n+1}$.

Let $\theta_0 > 0$ and assume $\ep_0 \in (0,1)$. 
Fix an arbitrary $\zo \in \hLa$. Replacing the group $\Gamma$ by a conjugate of its, 
we may assume that $\zo = (\xo, -e_1)$, where $\xo_1 = 1$ and $e_1 = (1,0,\ldots,0)\in \H^{n+1}$. 
Then $W^u_{\ep_0}(\zo)$ is a subset of 
$$W = \{ (x,-e_1)\in \sh : x_1 = 1\}\;,$$
and $E^u(\zo)$ can be naturally identified with $\pH = \{0\} \times \R^n$. 

In what follows for any $x= (x_1, \ldots,x_{n+1})\in \H^{n+1}$ we denote
$x' = (0,x_2, \ldots, x_{n+1}) \in \pH$.

Consider an arbitrary  $\hz = (\hx,-e_1) \in \hLa\cap W$ close to $\zo$, and
let $b\in \pH$ be a direction of $\hLa$-density at $\hz$, $\| b\| = 1$.  Setting
\be
\ep = \min \{ \ep_0/2, \theta_0/14\} < \frac{1}{2} \;,
\ee
the above implies the existence of
$(\hy, -e_1) \in \hLa\cap W\setminus \{ \hz\}$ such that
$$\|\hy - \hx\| < \ep\quad, \quad \left\|\frac{\hy-\hx}{\|\hy - \hx\|} - b \right\| < \ep\;.$$
Similarly, there exists $\tz = (\tx, -e_1) \in \hLa\cap W\setminus \{ \hz\}$ such that
$$\|\tx - \hx\| < \ep\| \hy-\hx\| \quad, \quad \left\|\frac{\tx-\hx}{\|\tx - \hx\|} - b \right\| 
< \ep\, \|\hy - \hx\|\;,$$
and $\tx'$ is a fixed point of a hyperbolic (loxodromic) element $\tg$ of $\Gamma$ 
(see e.g. Sect. 12.1 in 
\cite{kn:Ratc}). {\bf Fix $\hy$ and $\tz = (\tx, -e_1)$ with the above properties}. We then have
\be
\left\|\frac{\hy-\tx}{\|\hy - \tx\|} - b \right\| < 3\ep\;.
\ee

Next, changing the coordinate system in $\H^{n+1}$ if necessary we will assume that 
$\tx = (1,0,\ldots,0)$. Then $\tg$ has the form $\tg = k\, A$ for some $k > 0$, $k \neq 1$, 
and an orthogonal transformation $A$ in $\pH$. Replacing $\tg$ by $\tg^{-1}$ if necessary,
we will assume that $k > 1$. Considering the minimal $A$-invariant linear subspace of
$\pH = \R^n$ containing $b$, one derives that there exists an infinite  sequence 
$1 \leq m_1 < m_2 < \ldots < m_p < \ldots$ of integers such that $A^{m_p} b \to b$
as $p \to \infty$. Choose $p$ sufficiently large so that $\| A^{m_p}b - b\| < \ep$ and
set $m = m_p$ and $q = \tg^{m} (\hy')$. Since $\hy'\in L(\Gamma)$ and $L(\Gamma)$ is 
$\Gamma$-invariant,
we have $q\in L(\Gamma)$. We will assume $m = m_p$ is chosen so large that
$$\|q\| = k^{m}\|\hy'\| > 1\;.$$
With this choice of $q$ we have
\be
\| q/\|q\| - b\| < 4\ep\;.
\ee
Indeed, using the choice of $m$ and (7.2) it follows that
$$\|q/\|q\|  - b\|  \leq \|A^{m}\hy'/ \|\hy'\| - A^{m} b\| + \| A^{m}b - b\|
< \| \hy'/\| \hy'\| - b\| + \ep < 4\ep\;,$$
which proves (7.3). Fix $m$ and $q$ with the above properties.

Next, denote by $\tS$ the horosphere in $\H^{n+1}$ of radius $1/2$ at $0$ and  by $S_0$ the 
horosphere at $q$ externally tangent to $\tS$. Then $\tx \in \tS$ and $W^s_{\ep_0}(\tz)$ coincides 
with (a certain part of) the inward unit (with respect to the Poincar\'e metric) normal field to 
$\tS$. Let $R$ be the radius of $S_0$
and $u$ be the tangent point of $\tS$ and $S_0$. Then $\ty = (u,\xi) \in W^s_{\ep_0}(\tz)$
for some vector $\xi$ (assuming that $\|q\|$, and therefore $R$ is chosen sufficiently large),
and $W^u_{\ep_0}(\ty)$ coincides locally with the outward unit normal field to $S_0$.
Notice that $\lim_{t\to \infty} \phi_t(\ty) = 0\in L(\Gamma)$ and
$\lim_{t\to -\infty} \phi_t(\ty) = q\in L(\Gamma)$, so the definition of $\hLa$ implies
$\ty \in \hLa$.

Set $\ep'= \ep$ and consider an arbitrary $z = (x,-e_1)\in W^u_{\ep}(\tz)$; then $\|x'\| < \ep$.
Let $a\in \pH$ and $h\in \R$ be such that $\|a\| = 1$, $\la a,b\ra \geq \theta_0$ and $|h| < \ep$.
We will now show that (9.1) (see Sect. 9 below) holds with $\delta = \frac{\theta_0}{4\|q\|}$, $\ty_1 = \ty$ and $\ty_2 = \tz$. 
(Then $\Delta(\exp^u_z(v), \pi_{\ty_2}(z)) = \Delta(\exp^u_z(v), \pi_{\tz}(z)) = 0$.)

Let $S$ be the horosphere of radius $1/2$ at $x'$ ; then locally $W^s_{\ep_0}(z)$ coincides with
the inward unit normal field to $S$ (see Figure 1). So, for
$\sigma = \pi_{\ty}(z) = [z,\ty] = W^s_{\ep_0}(z) \cap \phi_{[-\ep,\ep]}(W^u_{\ep_0}(\ty))$
we have $\sigma = (v,\eta)$ for some $v \in S$ and 
$\phi_{t_1}(\sigma) \in W^u_{\ep_0}(\ty)$.  Thus, if $S_1$ is the horosphere at $x'$ tangent to $S_0$
(necessarily at the foot point of $\phi_{t_1}(\sigma) $) and $r_1$ is the radius of $S_1$,
then $t_1 = \ln(2r_1)$.  On the other hand, by elementary geometry,
$(R+r_1)^2 = \|q-x'\|^2 + (R-r_1)^2\;,$
so $r_1 = \|q-x'\|^2/(4R)$ and therefore
$\Delta(z,\ty) = t_1 = \ln \frac{\|q-x'\|^2}{2R}\;.$

In the same way for $\omega = \exp_z(ha) = (x+ha, -e_1)$ one obtains
$\Delta(\omega, \ty) = \ln \frac{\|q-x' - h a\|^2}{2R}\;.$
Therefore
$$\ttt = \Delta(\exp_z(ha),\pi_{\ty}(z)) = \Delta(\omega, \pi_{\ty}(z)) = 
\Delta(\omega, \ty) - \Delta(z, \ty)  = \ln  \frac{\|q-x' - h a\|^2}{\|q-x'\|^2}\;.$$
Using the fact that $|\ln (1+x)| \geq |x|/2$ for $|x| < 1$, one gets
$$|\ttt| = \left| \ln\left[ 1- \frac{2h}{\|q-x'\|}\, 
\left\langle \frac{q-x'}{\|q-x'\|}, a\right\rangle + \frac{h^2}{\|q-x'\|^2} \right] \right|
\geq \frac{|h|}{2\|q-x'\|} \left| 2 \left\langle \frac{q-x'}{\|q-x'\|}, a\right\rangle 
-\frac{h}{\|q-x'\|} \right|\;.$$
Now (7.3)   implies
$$\left\| \frac{q-x'}{\|q-x'\|} - b\right\| < 4\ep + \left\|\frac{q-x'}{\|q-x'\|} 
- \frac{q}{\|q\|}\right\| \leq 4\ep + 2 \frac{\|x'\|}{\|q\|} < 6\ep\;,$$
and using $\la a,b\ra \geq \theta_0$ we get
$$2 \left\langle \frac{q-x'}{\|q-x'\|}, a\right\rangle = 
2\la b,a\ra + 2 \left\langle \frac{q-x'}{\|q-x'\|} -b, a\right\rangle \geq 2\theta_0 -12 \ep\;.$$
Moreover, $\|q-x'\| \geq 1-\ep > 1/2$, so  $|h|/\| q-x'\| < 2\ep$ which combined with the above,
$\|q-x'\| \leq \|q\|+\ep < 2\|q\|$  and (7.1) gives
$$|\ttt| \geq \frac{|h|}{4 \|q\|} \, (2\theta_0 -14\ep) \geq \frac{\theta_0}{4 \|q\|}\, 
|h| = \delta\, |h|$$
for all $h$ with $|h| < \ep$. 

This proves that (LNIC) is fulfilled, thus completing the proof of Lemma 3.

\section{Open billiard flows}
\renewcommand{\theequation}{\arabic{section}.\arabic{equation}}
\setcounter{equation}{0}

In this section we prove Theorem 4. Let $K$ be a subset of ${\R}^{N}$ ($N\geq 2$) of the form
$K = K_1 \cup K_2 \cup \ldots \cup K_{\kappa_0},$where $K_i$ are compact 
strictly convex disjoint domains in $\R^{N}$ with 
$C^r$ ($r \geq 3$) {\it boundaries} $\Gamma_i = \dk_i$ and $\kappa_0 \geq 3$. 
Set $\Omega = \overline{{\R}^N \setminus K}.$ 
Throughout this section we assume that $K$ satisfies the following (no-eclipse) condition: 
$${\rm (H)} \quad \quad\qquad  
\begin{cases}
\mbox{\rm for every pair $K_i$, $K_j$ of different connected components 
of $K$ the convex hull }\cr
\mbox{\rm of $K_i\cup K_j$ has no
common points with any other connected component of $K$. }\cr
\end{cases}$$
With this condition, the {\it billiard flow} $\varphi_t$ defined on the {\it cosphere bundle} $S^*(\Omega)$ 
in the standard way is called an open billiard flow.
It has singularities, however its restriction to the {\it non-wandering set} $\Lambda$ has only 
simple discontinuities at reflection points.  
Moreover, $\Lambda$  is compact, $\varphi_t$ is hyperbolic and transitive
on $\Lambda$, and  it follows from  \cite{kn:St1} that $\varphi_t$ is  non-lattice and therefore by  
a result of  Bowen \cite{kn:B}, it is topologically weak-mixing on $\Lambda$.

Denote by $A$ the $\kappa_0\times \kappa_0$ matrix with entries $A(i,j) = 1$ if $i \neq j$ and
$A(i,i)= 0$ for all $i$, and define $\sa$ and $\saa$ as in Sect. 1. Given $\xi \in \sa$, let
$\ldots, P_{-2}(\xi), P_{-1}(\xi), P_0(\xi), P_1(\xi), P_2(\xi) , \ldots$
be the successive reflection points of the unique billiard trajectory in the exterior of $K$ such that
$P_j(\xi) \in K_{\xi_j}$ for all $j \in \Z$. Set
$f(\xi) = \| P_0(\xi) - P_1(\xi)\|\;,$ and define the map
$\Phi : \sa \longrightarrow \mtb = \mt \cap S^*_\mt(\Omega)$ by 
$$\Phi(\xi) = (P_0(\xi), (P_1(\xi) - P_0(\xi))/ \|P_1(\xi) - P_0(\xi)\|)\;.$$
Then $\Phi$ is a bijection such that $\Phi\circ \sigma = B \circ \Phi$, where 
$B : \mtb \longrightarrow \mtb$ is the {\it billiard ball map}.
Choosing appropriately $\theta \in (0,1)$, we have
$f  \in {\mathcal F}_{\theta^2}(\sa)$  (see e.g. \cite{kn:I}).

By Sinai's Lemma (see e.g. \cite{kn:PP}), there exists a function $\tf \in \ff_\theta(\sa)$ 
depending on future coordinates only and $\chi_f \in \ff_\theta(\sa)$ such that
$f(\xi) = \tf(\xi) + \chi_f(\xi) - \chi_f(\sigma \xi)$ for all $\xi \in \sa$.
As in  the proof of Sinai's Lemma, for any $k = 1, \ldots, k_0$ choose and fix an arbitrary sequence 
$\eta^{(k)} = (\ldots, \eta^{(k)}_{-m}, \ldots, \eta^{(k)}_{-1}, \eta^{(k)}_0)\in \Sigma_a^-$ 
with $\eta^{(k)}_0 = k$. Then for any $\xi \in \sa$ (or $\xi \in \saa$) set
$$e(\xi) = (\ldots, \eta^{(\xi_0)}_{-m}, \ldots, \eta^{(\xi_0)}_{-1},
\eta^{(\xi_0)}_0= \xi_0, \xi_1, \ldots, \xi_m, \ldots)\in \sa \;.$$ 
Then we have
$$\chi_f(\xi) = \sum_{n=0}^\infty [ f(\sigma^n (\xi)) - f(\sigma^n e(\xi))]\;.$$
As before, let $P = P_f \in \R$ be such that $\Pr(-P\, \tf) = 0$ (then $\Pr(-P\, f) = 0$
as well), let $m_0$ be the equilibrium state of $-P_f \tf$, and let $\alpha = \int_{\saa} f\, dm_0$.

Next, let $\rr = \{ R_i\}_{i=1}^k$ be a Markov family of rectangles with $R_i = [U_i  , S_i ]$,
$U_i \subset W^u_\ep(z_i) \cap \mt$ and $S_i \subset W^s_\ep(z_i)\cap \mt$, respectively, 
for some $z_i\in \mt$ (see Sect. 6 above).  
Taking $\chi$ sufficiently small, we may assume that each rectangle $R_i$ is `between
two boundary components' $\Gamma_{p_i}$ and $\Gamma_{q_i}$ of $K$, that is for any
$x\in R_i$, the first backward reflection point of the billiard trajectory $\gamma$ determined by $x$
belongs to $\Gamma_{p_i}$, while the first forward reflection point of $\gamma$ belongs to $\Gamma_{q_i}$.
Moreover, using the fact that the intersection of $\mt$ with each
cross-section to the flow $\varphi_t$ is a Cantor set, we may assume that the Markov family
$\rr$ is chosen in such a way that, apart from the standard properties (a) and (b) in Sect. 6, it also
satisfies the following:

\ms

(c)  for any $i = 1, \ldots, k$ we have $\partial_\mt U_i  = \e$.

\ms

Finally, partitioning every $R_i$ into finitely many smaller rectangles, cutting $R_i$ along some 
unstable leaves, and removing some rectangles from the family formed in this way, we may assume that

\ms

(d) for every $x\in R$ the billiard trajectory of $x$ from $x$ to $\pp(x)$ makes exactly one
reflection.

\ms

\def\piU{\pi^{(U)}}

From now on we will assume that $\rr = \{ R_i\}_{i=1}^k$ is a fixed Markov family for  $\varphi_t$
of size $\chi < \ep_0/2$ satisfying the conditions (a), (b) from Sect. 6 and the above conditions
(c) and (d). Define $U = \cup_{i=1}^k U_i$ and $\sigma : U   \longrightarrow U$ as in Sect. 6. 
As in Sect. 6, we will assume that $\rr$ is chosen so that $\tau$ is non-lattice.

Under the conditions in Theorem 4 it follows from Theorem 6 in Sect. 9 below that
the Ruelle transfer operators related to $\tau$ are weakly contracting. 
This allows  to apply Theorem 6. One particular case when these conditions are satisfied concerns the following 
{\it pinching condition}:

\ms

\noindent
{\sc (P)}:  {\it There exist  constants $C > 0$ and $\alpha > 0$ such that for every $x\in \mt$
we have
$$\frac{1}{C} \, e^{\alpha_x \,t}\, \|u\| \leq \| d\varphi_{t}(x)\cdot u\| 
\leq C\, e^{\beta_x\,t}\, \|u\| \quad, \quad  u\in E^u(x) \:\:, t > 0 \;,$$
for some constants $\alpha_x, \beta_x > 0$ depending on $x$ but independent of $u$ with
$\alpha \leq \alpha_x \leq \beta_x $ and $2\alpha_x - \beta_x \geq \alpha$ for all $x\in \mt$.}

\ms

Notice that when $N = 2$ this condition is always satisfied. For $N \geq 3$, (P) follows from 
certain estimates on the eccentricity of the connected components $K_j$ of $K$ -- see 
\cite{kn:St4} for a more precise result. It turns out that for $n \geq 3$ the condition (P) 
is always satisfied when the minimal distance between  distinct connected 
components of $K$ is relatively large compared to the maximal sectional curvature of $\partial K$.
According to general regularity results (\cite{kn:H}), (P) implies that $W^u_\ep(x)$ and $W^s_\ep(x)$
are $C^{1+\delta}$ in $x\in \mt$ for some $\delta > 0$. This and the main result in
\cite{kn:St4} imply the following

\begin{prop}
Assume that the billiard flow $\varphi_t$ satisfies the pinching condition (P) on $\mt$. 
Then the Ruelle transfer operators related to $\tf$ are weakly contracting.
\end{prop}

\ms

\noindent
{\it Proof of Theorem} 4. Assume that the conditions of Theorem 4 are satisfied.
Then, as mentioned above, $\tau$ is non-lattice and the Ruelle transfer operators 
related to $\tau$ are weakly contracting, so we can apply Theorem 5 from Sect. 6.

Let $\aa$ be the matrix defined in Sect. 6 using the Markov family $\rr$. As in Sect. 2
in \cite{kn:PeS} one defines a natural bijection $\ss : \sAA \longrightarrow \saa$
which commutes with the shifts. Apart from that there is a natural map 
$\W: U \longrightarrow \sAA$  such that
$\sigma \circ \W = \W \circ \sigma$. Let $P_\tau\in \R$ be such that
$\Pr(-P_\tau\, \tau) = 0$. It is easy to see that $P_\tau = P_f$. Indeed, 
first notice that the map $\W$ is continuous (and therefore a homeomorphism) 
when $U$ is considered with the Riemannian metric and $\sAA$ with the metric $d_\theta$,
so $\Pr (-P_\tau\, \tau\circ \W^{-1}) = 0$ (see e.g. Theorem 9.8 in \cite{kn:Wal}). Next, 
for $r = \tau\circ \W^{-1}$ it follows from 
(3.4) in Sect.3 in \cite{kn:PeS} that there exists a continuous function
$\mu : \saa \longrightarrow \R$ such that 
$r  = \tf \circ \ss + \mu \circ \sigma - \mu$.
Thus, $\Pr (a\, r) = \Pr(a\, \tf)$ for any $a\in \R$ (see e.g. \cite{kn:PP}),
so in particular, $\Pr(-P_\tau\, \tf) = 0$, and therefore  $P_\tau = P_f$.
In a similar way we see that if $m'_0$ is the equilibrium state of $- P_\tau \,\tau$
on $U$, then $\alpha = \int_{\saa} \tf\, dm_0 = \int_U \tau\, dm'_0$. 

It remains to notice that if $\sigma^n(x) = x$ for some $x \in U$, and if $n$ is the smallest  integer with this property, then $x$ generates a periodic billiard orbit $\gamma$ with $n$ reflection points
and  $T_\gamma = \tau^n(x)$. Every periodic billiard orbit with $n$ reflection points
is obtained in this way, and we get the same orbit from $n$ different $x$.
With this remark, using Theorem 5 from Sect. 6, we get the 
estimates (\ref{eq:1.11}) and (\ref{eq:1.12}).
\endofproof

\section{Spectral estimates for Ruelle transfer operators}
\renewcommand{\theequation}{\arabic{section}.\arabic{equation}}
\setcounter{equation}{0}

Let again $\varphi_t : M \longrightarrow M$ be a $C^2$ Axiom A flow and $\mt$
be  a basic set for $\phi_t$.
For any $x \in \mt$, $T > 0$ and $\delta\in (0,\ep]$ set
$$B^u_T (x,\delta) = \{ y\in W^u_{\ep}(x) : d(\varphi_t(x), \varphi_t(y)) \leq \delta \: \: 
, \:\:  0 \leq t \leq T \}\;.$$

We will say that $\varphi_t$ has a {\it regular distortion along unstable manifolds} over
the basic set $\mt$  if there exists a constant $\ep_0 > 0$ with the following properties:

\ms 

(a) For any  $0 < \delta \leq   \ep \leq \ep_0$ there exists a constant 
$R =  R (\delta , \ep) > 0$ such that 
$$\diam( \mt \cap B^u_T(z ,\ep))   \leq R \, \diam( \mt \cap B^u_T (z , \delta))$$
for any $z \in \mt$ and any $T > 0$.

\ms

(b) For any $\ep \in (0,\ep_0]$ and any $\rho \in (0,1)$ there exists $\delta  \in (0,\ep]$
such that for  any $z\in \mt$ and any $T > 0$ we have
$\diam ( \mt \cap B^u_T(z ,\delta))   \leq \rho \; \diam( \mt \cap B^u_T (z , \ep))\;.$

\bs

Part (a) of the above condition resembles the Second Volume Lemma of Bowen and Ruelle 
\cite{kn:BR} about balls in 
Bowen's metric; this time however we deal with diameters instead of volumes.
Sect. 8 in \cite{kn:St3} describes a rather general class of flows on basic sets satisfying this condition.
In fact, there are reasons to believe that this may actually hold for all $C^2$
flows on basic sets -- see the comments in Sect. 1 in \cite{kn:St3}.

In the special case when the flow satisfies the pinching condition (P) over $\mt$ (see Sect. 8, where it is stated for open billiard flows; for general flows on basic sets it is similar), it
follows from Theorem 7.1 in \cite{kn:St3} that
$\varphi_t$ has a regular distortion along unstable manifolds over $\mt$. As we mentioned in
Sect. 8 above, when the local unstable manifolds are one-dimensional (P) is 
always satisfied. For open billiards (see Sect. 8 again) the condition (P) is always satisfied 
when the minimal distance between  distinct connected components of $K$ is relatively large 
compared to the maximal sectional curvature of $\dk$. 
An analogue of the latter for manifolds $M$ of strictly negative curvature would be to require that 
the sectional curvature is between $-K_0$ and $-a\, K_0$ for some constants $K_0 > 0$ and 
$a\in (0,1)$.  It follows from the arguments in  \cite{kn:HP} that when 
$a = 1/4$  the geodesic flow on $M$ satisfies the pinching condition (P).

In what follows we deal with flows $\varphi_t$ over basic sets $\mt$ having a regular distortion 
along unstable manifolds. Apart from that, we impose an additional
{\it local non-integrability condition} (LNIC) which we state below.

It follows from the hyperbolicity of $\mt$  that if  $\epsilon_0 > 0$ is sufficiently small,
there exists $\ep_1 > 0$ such that if $x,y\in \mt$ and $d (x,y) < \ep_1$, 
then $W^s_{\ep_0}(x)$ and $\varphi_{[-\ep_0,\ep_0]}(W^u_{\ep_0}(y))$ intersect at exactly 
one point $[x,y ] \in \mt$  (cf. \cite{kn:KH}). That is, there exists a unique 
$t\in [-\ep_0, \ep_0]$ such that $\varphi_t([x,y]) \in W^u_{\ep_0}(y)$. Setting 
$\Delta(x,y) = t$, defines the so called {\it temporal distance
function}. For $x, y\in \mt$ with $d (x,y) < \ep_1$, define
$\pi_y(x) = [x,y] = W^s_{\ep}(x) \cap \varphi_{[-\ep_0,\ep_0]} (W^u_{\ep_0}(y))\;.$
Thus, for a fixed $y \in \mt$, $\pi_y : W \longrightarrow \varphi_{[-\ep_0,\ep_0]} (W^u_{\ep_0}(y))$ is the
{\it projection} along local stable manifolds defined on a small open neighborhood $W$ of $y$ in $\mt$.

Given $z \in \mt$, let $\exp^u_z : E^u(z;\ep_0) \longrightarrow W^u_{\ep_0}(z)$  
be the corresponding {\it exponential map}.
A  vector $b\in E^u(z)\setminus \{ 0\}$ will be called  {\it tangent to $\mt$} at
$z$ if there exist infinite sequences $\{ v^{(m)}\} \subset  E^u(z)$ and 
$\{ t_m\}\subset \R\setminus \{0\}$
such that $\exp^u_z(t_m\, v^{(m)}) \in \mt \cap W^u_{\ep}(z)$ for all $m$, $v^{(m)} \to b$ and 
$t_m \to 0$ as $m \to \infty$. 
It is easy to see that a vector $b\in E^u(z)\setminus \{ 0\}$ is  tangent to $\mt$ at
$z$ if there exists a $C^1$ curve $z(t)$ ($0\leq t \leq a$) in $W^u_{\ep}(z)$ for some $a > 0$ 
with $z(0) = z$ and $\dot{z}(0) = b$ such that $z(t) \in \mt$ for arbitrarily small $t >0$.

The following is the  {\it local non-integrability condition} for $\varphi_t$ and $\mt$ mentioned
above.

\medskip

\noindent
{\sc (LNIC):}  {\it There exist $z_0\in \mt$,  $\ep_0 > 0$ and $\theta_0 > 0$ such that
for any  $\ep \in (0,\ep_0]$, any $\hz\in \mt \cap W^u_{\ep}(z_0)$  and any tangent vector 
$\eta \in E^u(\hz)$ to $\mt$ at $\hz$ with 
$\|\eta\| = 1$ there exist  $\tz \in \mt \cap W^u_{\ep}(\hz)$, 
$\ty_1, \ty_2 \in \mt \cap W^s_{\ep}(\tz)$ with $\ty_1 \neq \ty_2$,
$\delta = \delta(\tz,\ty_1, \ty_2) > 0$ and $\ep'= \ep'(\tz,\ty_1,\ty_2)  \in (0,\ep]$ such that
\be
|\Delta( \exp^u_{z}(v), \pi_{\ty_1}(z)) -  \Delta( \exp^u_{z}(v), \pi_{\ty_2}(z))| \geq \delta\,  \|v\| 
\ee
for all $z\in W^u_{\ep'}(\tz)\cap\mt$  and  $v\in E^u(z; \ep')$ with  $\exp^u_z(v) \in \mt$ and
$\la \frac{v}{\|v\|} , \eta_z\ra \geq \theta_0$,   where $\eta_z$ is the parallel 
translate of $\eta$ along the geodesic in $W^u_{\ep_0}(z_0)$ from $\hz$ to $z$.}

\ms

One would  expect that (LNIC) is satisfied in most interesting cases. For example, it was shown in 
\cite{kn:St4} that open billiard flows (in any dimension) with $C^1$ (un)stable laminations over the 
non-wandering set $\mt$ always satisfy (LNIC).

If $\varphi_t$ is a $C^2$ contact flow on $M$, i.e. there exists a $C^2$ invariant 
one-form $\omega$ such that $\omega \wedge (d\omega)^n$ is a volume form on $M$, where
$\dim(M) = 2n+1$, then the following condition (ND) implies (LNIC) (see Proposition 6.1
in \cite{kn:St3}).

\ms

\noindent
{\sc (ND)}:  {\it There exist $z_0\in \mt$, $\ep > 0$  and  $\mu_0 > 0$ such that
for any $\ep \in (0,\ep_0]$, any $\hz \in \mt \cap W^u_{\ep}(z_0)$ and any unit vector 
$\eta \in E^u(\hz)$ tangent to 
$\mt$ at $\hz$ there exist $\tz \in \mt \cap W^u_{\ep}(\hz)$, $\ty \in W^s_\ep(\tz)$ and   
a unit vector $\xi \in E^s(\ty)$ tangent to $\mt$ at $\ty$ with 
$|d\omega_{\tz}(\xi_{\tz},\eta_{\tz}) | \geq \mu_0$,
where $\eta_{\tz}$ is the parallel translate of $\eta$ along the geodesic in $W^u_{\ep}(\tz)$ from
$\hz$ to $\tz$, while $\xi_{\tz}$ is the parallel translate of $\xi$ along the geodesic in 
$W^s_{\ep}(\tz)$ from $\ty$ to $\tz$.}

\ms

\noindent
{\bf Remark.} It appears the above condition would become significantly more
restrictive if one requires the existence of a unit vector $\xi \in E^s(\tz)$ tangent to $\mt$ at 
$\tz$ with  $|d\omega_{\tz}(\xi,\eta_{\tz}) | \geq \mu_0$. The reason for this is that in general the set
of unit tangent vectors to $\mt$ does not have to be closed in the bundle $E^s_\mt$.
That is, there may exist a point  $\tz \in \mt$, a
sequence $\{ z_m\} \subset W^s_{\ep}(\tz)\cap \mt$ and for each $m$ a unit vector
$\xi_m$ tangent to $\mt$ at $z_m$ such that $z_m \to z$ and $\xi_m \to \xi$ as $m \to \infty$,
however $\xi$ is not tangent to $\mt$ at $\tz$. A similar comment can be made about (LNIC), where
requiring $\ty_2 = \tz$ would replace (9.1) by $|\Delta( \exp^u_{z}(v), \pi_{\ty}(z))| \geq \delta\,  \|v\| $
with $\ty = \ty_1$, which is still a rather general non-integrability condition. However in its present form
(LNIC) is a substantially weaker condition.

\ms

Given a Lipschitz real-valued function $f$  on $U$, set $g = g_f = f - P\tau$, where 
$P = P_f\in \R$ is the unique 
number such that the topological pressure $\Pr_\sigma(g)$ of $g$ with respect to $\sigma$ is 
zero (cf. e.g. \cite{kn:PP}). For $a, b\in \R$, consider the {\it Ruelle transfer operator}
$L_{g-(a+\i b)\tau}$ on the space $\clip (U)$ of Lipschitz functions $g: U \longrightarrow \C$. 
By  $\Lip(g)$ we denote the Lipschitz constant of $g$ and  by $\| g\|_\infty$ the 
{\it standard $\sup$ norm}   of $g$ on $U$. As in Sect. 6 above, we will use the norm $\|.\|_{{\lip},b}$ on $\clip (U)$ defined by 
$\| h\|_{{\lip},b} = \|h\|_\infty + \frac{\Lip(h)}{|b|}$.

The following result from \cite{kn:St3} has been used several times in  previous sections.

\begin{thm}
$($\cite{kn:St3}$)$  Let $\varphi_t : M \longrightarrow M$ be  a $C^2$ Axiom A flow on
a $C^2$ complete  Riemann manifold satisfying the condition (LNIC) and having a
regular distortion along unstable manifolds over a  basic set $\mt$. Assume in addition  that the local 
holonomy maps along stable laminations through $\mt$ are uniformly Lipschitz. Then for any  Lipschitz 
real-valued function $F$ we have the following: for every $\epsilon > 0$ there exist constants 
$0 < \rho < 1$, $a_0 > 0$ and  $C > 0$ such that if $a,b\in \R$  satisfy $|a| \leq a_0$ and $|b| \geq 1/a_0$, then 
$$\|\lc_{F -(P(F)+ a + \i b)\tau}^m h \|_{{\lip},b} \leq C \;\rho^m \;|b|^{\ep}\; \| h\|_{{\lip},b}\;$$
for every integer $m > 0$ and every  $h\in \clip (U)$. 
In particular the spectral radius  of $\lc_{F-(P(F)+ a +\i b)\tau}$ on $\clip(U)$ does not exceed   $\rho$.
\end{thm}

As an immediate consequence of this theorem we get the following (see \cite{kn:D} or Corollary 3.3 in \cite{kn:St2}):

\begin{cor}
Under the assumptions of Theorem 6, the Ruelle transfer operators related to $\tau$ are weakly contracting.
\end{cor}

\includegraphics[scale= 1]{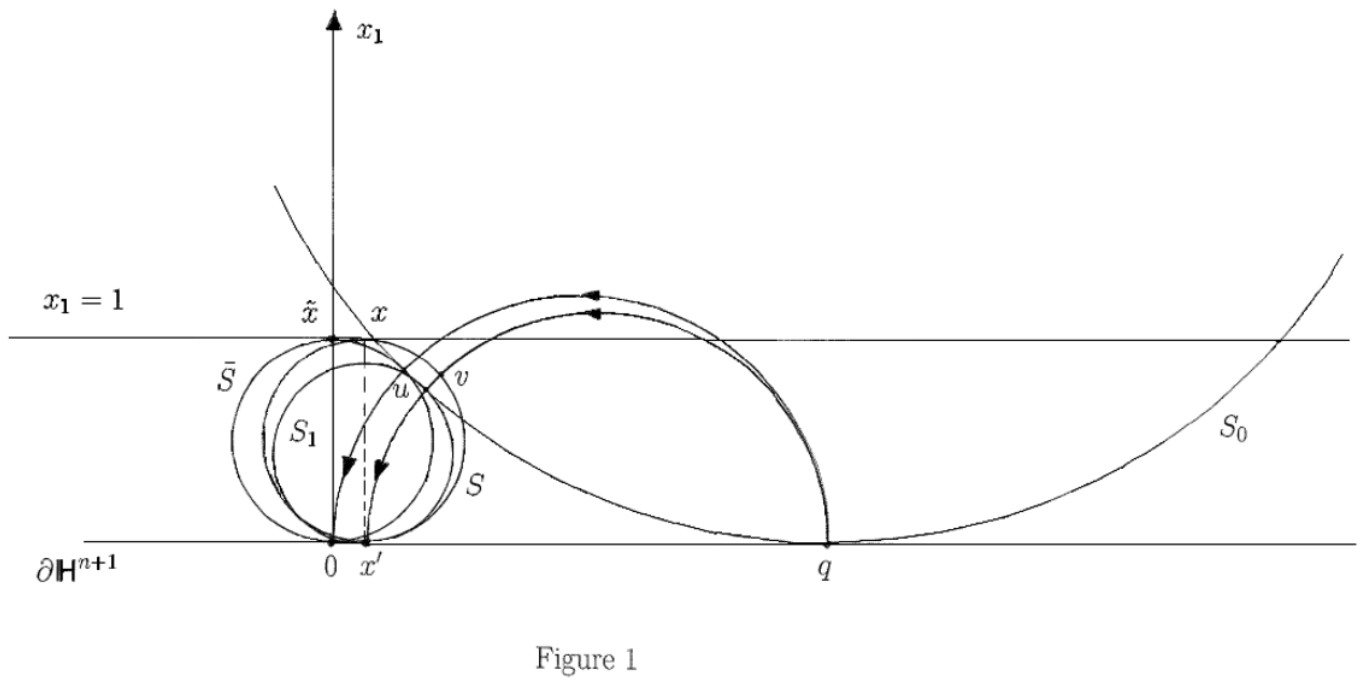}

\end{document}